\documentclass[12pt]{article}
\usepackage{theorem,exscale}
\usepackage{amssymb}
\usepackage[lined,algoruled,linesnumbered]{algorithm2e}

\usepackage{latexsym}
\usepackage{fullpage}
\usepackage{float}
\usepackage{algorithmic}
\usepackage{epsf}
\usepackage{subeqn}
\usepackage{multicol}

\usepackage{hyperref}

\def\R{\mathbb{R}}

\def\eqref#1{{\normalfont(\ref{#1})}}
\def\SDPt{{\bf SDP\,}}
\def\SDP{\mbox{\boldmath$SDP\,$}}
\def\SDPc{{\mbox{\boldmath$SDP,\,$}}}
\def\EDM{\mbox{\boldmath$EDM\,$}}
\def\SNL{\mbox{\boldmath$SNL\,$}}
\def\SNLt{{\bf SNL\,}}
\def\SNLc{\mbox{\boldmath$SNL,\,$}}
\def\SNLLS{\mbox{\boldmath$SNL_{LS}\,$}}
\def\SNLM{\mbox{\boldmath$SNL_{M}\,$}}
\def\SNLMV{\mbox{\boldmath$SNL_{MV}\,$}}

\def\EDMt{\bf EDM\,}

\def\EDMC{\mbox{\boldmath$EDMC\,$}}

\def\SNLMN{\mbox{\boldmath$SNL_{MN}\,$}}
\def\SNLMND{\mbox{\boldmath$SNL_{MV}-D\,$}}

\newtheorem{lem}{Lemma}[section]
\newtheorem{assumption}{Assumption}[section]

\newtheorem{thm}{Theorem}[section]
\newtheorem{cor}{Corollary}[section]
\newtheorem{rem}{Remark}[section]

\newcommand{\diag}{{\rm diag\,}}
\newcommand{\offDiag}{{\rm offDiag\,}}
\newcommand{\Diag}{{\rm Diag\,}}
\newcommand{\Ss}{{\mathcal S}}
\newcommand{\Sn}{{\mathcal S}^n }
\newcommand{\Snm}{{\mathcal S}^{n-1} }

\newcommand{\KK}{{\mathcal K} }
\newcommand{\LL}{{\mathcal L} }
\newcommand{\NN}{{\mathcal N} }
\newcommand{\MM}{{\mathcal M} }

\newcommand{\ZZ}{{\mathcal Z} }
\newcommand{\YY}{{\mathcal Y} }
\newcommand{\VV}{{\mathcal V} }
\newcommand{\EE}{{\mathcal E} }
\newcommand{\FF}{{\mathcal F} }
\newcommand{\DD}{{\mathcal D} }

\newcommand{\GG}{{\mathcal G} }

\newcommand{\RR}{{\mathcal R} }

\newcommand{\trace}{{\rm trace\,}}

\newcommand{\relint}{{\rm relint\,}}
\newcommand{\argmin}{{\rm argmin\,}}
\newcommand{\rank}{{\rm rank\,}}

\newcommand{\tr}{{\rm trace\,}}
\newcommand{\svec}{{\rm svec\,}}
\newcommand{\sMat}{{\rm sMat\,}}
\newcommand{\sblk}{{\rm sblk\,}}
\newcommand{\sBlk}{{\rm sBlk\,}}

\newcommand{\kvec}{{\rm vec\,}}
\newcommand{\Mat}{{\rm Mat\,}}
\newcommand{\bpr}{{\bf Proof.} \hspace{1 em}}
\newcommand{\QED}{\hfill ~\rule[-1pt] {8pt}{8pt}\par\medskip ~~}
\newcommand{\epr}{\QED}

\newcounter{count}
\newcommand{\beqr}{\addtocounter{count}{1} \begin{eqnarray}}

\renewcommand{\theequation}{\thesection.\thecount}
\newcommand{\beq}{\addtocounter{count}{1} \begin{equation}}
\newcommand{\bet}{\addtocounter{count}{1} \begin{table}}
\newcommand{\eeq}{ \end{equation} }
\newcommand{\eeqr}{ \end{eqnarray} }

\newcommand{\bt}{ \begin{tabular} }
\newcommand{\et}{ \end{tabular} }

\newcommand{\ip}[2]{\left\langle #1, #2 \right\rangle}

\begin{document}

\bibliographystyle{plain}
\title{
Sensor Network Localization, Euclidean Distance Matrix Completions,
 and Graph Realization
}
             \author{
Yichuan Ding
   \thanks{Research supported by Natural Sciences Engineering Research
     Council Canada. E-mail y7ding@math.uwaterloo.ca}
\and
Nathan Krislock
   \thanks{Research supported by Natural Sciences Engineering Research
     Council Canada. E-mail ngbkrislock@uwaterloo.ca}
\and
Jiawei Qian
   \thanks{Research supported by Natural Sciences Engineering Research
     Council Canada. E-mail j2qian@student.cs.uwaterloo.ca}
\and
Henry Wolkowicz
   \thanks{Research supported by Natural Sciences Engineering Research
     Council Canada. E-mail hwolkowicz@uwaterloo.ca}
}
\date{\today}
          \maketitle
\begin{center}
          University of Waterloo\\
          Department of Combinatorics and Optimization\\
          Waterloo, Ontario N2L 3G1, Canada\\
          Research Report CORR 2006-23
\end{center}

\vspace{0.1in}

\noindent {\bf AMS Subject Classification:}

\tableofcontents
\listoftables
\listoffigures

\begin{abstract}
\noindent
We study Semidefinite Programming, \SDPc relaxations for 
Sensor Network Localization, \SNLc with anchors and 
with noisy distance information.
The main point of the paper is to view \SNL as a 
(nearest) Euclidean Distance Matrix, \EDM, completion problem and to show the
advantages for using this latter, well studied model.
We first show that the current popular \SDP relaxation is equivalent to 
known relaxations in the literature for \EDM completions.
The existence of anchors in the problem is {\em not} special. The 
set of anchors simply
corresponds to a given fixed clique for the graph of the \EDM problem. 
We next propose a method of projection when a large clique 
or a dense subgraph is identified 
in the underlying graph. This projection reduces the size, and
improves the stability, of the relaxation.

In addition, viewing the problem as an \EDM completion problem yields
better low rank approximations for the low dimensional realizations.

And, the projection/reduction 
procedure can be repeated for other given cliques of sensors or for sets of
sensors, where many distances are known. Thus, further size reduction can
be obtained.

Optimality/duality conditions
and a primal-dual interior-exterior path following algorithm are
derived for the \SDP relaxations
We discuss the relative stability and strength of two
formulations and the corresponding algorithms that are used.
In particular, we show that the quadratic formulation arising from the
\SDP relaxation is better conditioned than the linearized form, that
is used in the literature and that arises from applying a Schur complement.
\end{abstract}

\section{Introduction}
\label{sect:intro}

We study ad hoc wireless sensor networks and the
sensor network localization, \SNL, problem with
anchors. The anchors have fixed known locations and the sensor-sensor and
sensor-anchor distances are known (approximately) 
if they are within a given (radio) range.
The problem is to approximate the positions of all the sensors, given that
we have only this partial information on the distances.
We use semidefinite programming, \SDP, relaxations to 
find approximate solutions to this problem.

In the last few years,
there has been an increased interest in the \SNL problem with anchors.
In particular, \SDP relaxations have been introduced that are specific
to the problem with anchors.
In this paper we emphasize that the existence of anchors is not special. 
The \SNL problem with anchors can be modelled as a 
(nearest) Euclidean Distance Matrix, 
\EDM, completion problem, a well studied problem.
There is no advantage to considering the anchors separately to other
sensors. The only property that distinguishes the anchors is that the 
corresponding set of nodes yields a clique in the graph.
This results in the failure of the Slater constraint qualification for the \SDP
relaxation.  
We then  show that we can take advantage of this liability. We can find
the smallest face of the \SDP cone that contains the feasible set and
project the problem onto this face.

This projection technique yields an equivalent
smaller dimensional problem, where the Slater constraint qualification holds.
Thus the problem size is reduced and the problem stability is improved.
In addition, viewing the problem as an \EDM completion leads to improved
low rank factorizations for the low dimensional realizations.
And,  by treating the anchors this way, we show that other
cliques of sensors or dense parts of the graph can similarly result in
a reduction in the size of the problem. In addition, not treating other cliques
this way can result in instability, due to loss of the Slater constraint
qualification.

We also derive optimality and duality conditions for the \SDP relaxations.
This leads to a primal-dual interior-exterior path following algorithm.
We discuss the robustness and stability of two
approaches. One approach is based on the quadratic constraint in matrix
variables that arises from the \SDP relaxation. 
The other approach uses the linearized version that is used in the
literature, and that is obtained from an application of the Schur
complement.  
Numerical tests comparing these two equivalent formulations
of the \SDP relaxation are included. They show that the quadratic
formulation is better conditioned and requires fewer iterations to reach
a desired relative duality gap tolerance.
These tests confirm results in the
literature, see \cite{GulerTuncel:98,ChuaTu:05},
on the conditioning of the central path and the comparison of different
barriers.

\subsection{Related Work and Applications}
The geometry of \EDM has been extensively studied in the literature, e.g.
\cite{MR86j:62133,MR90a:92082} and more recently in \cite{AlKaWo:97,homwolkA:04} 
and the references therein. The latter two references studied algorithms
based on \SDP formulations of the \EDM completion problem.

Several recent papers have developed algorithms for
the \SDP relaxation designed specifically for \SNL with anchors, e.g.
\cite{BiswasYe:04,Jin:05,MR2191577,biswasliangtohwangye,SoYe:05,biswasliangtohye:05,WangZhengBoydYe:06,KrPiWo:06}.
Relaxations using second order cones are studied in e.g.
\cite{tseng:04,Tsengsiam:05}.

The \SDP relaxations solve a closest \SDP matrix problem
and generally use the $\ell_1$ norm.
The $\ell_2$ norm is used in
\cite{KrPiWo:06}, where 
the noise in the radio signal is assumed to come
from a multivariate normal distribution with
mean $0$ and variance-covariance matrix $\sigma^2 I$,  i.e. from
a spherical normal distribution so that
the least squares estimates are the maximum likelihood estimates.
(We use the $\ell_2$ norm as well in this paper.
Our approach follows that in \cite{AlKaWo:97} for \EDM completion {\em
without} anchors.)

Various applications for \SNL are discussed in the references mentioned
above. These applications include e.g.
natural habitat monitoring, earthquake detection, and weather/current
monitoring.

\subsection{Outline}
The formulation of the \SNL problem as both a feasibility question 
and as a least squares approximation is presented in Section
\ref{sect:formulation}.
We continue in Section \ref{sect:distgeom} with background, notation,
including information on the linear transformations and adjoints used in the
model.  In particular, since this paper emphasizes using \EDM,
this section provides details on {\em distance geometry}. In particular,
we provide details on the linear mappings between \EDM and \SDP matrices.

The \SDP relaxations are presented in Section
\ref{sect:relax}. This section contains the details
for the four main contributions of the paper: i.e.
\begin{quote}
(i) the connection of \SNL with \EDM; (ii) the projection technique for
cliques and dense sets of sensors; (iii) the improved approximation scheme for 
locating the sensors from the \SDP relaxation; and (iv) a numerical
comparison showing the better conditioning of the quadratic formulation
relative to the linear formulation used in the literature.
\end{quote}
We begin in Section \ref{sect:relaxcurrent} with
several lemmas that describe the feasible set of the \SDP relaxation,
e.g. Lemma \ref{lem:equivcontr} provides several equivalent
characterizations that show the connection with \EDM.
The key to the connection is the loss of the Slater constraint
qualification (strict feasibility); but one can project onto the 
{\em minimal face} in order to
obtain the Slater condition and guarantee numerical 
stability and strong duality. This Lemma also shows the equivalent
representations of the feasible set by a quadratic and a linear
semidefinite constraint.
Then Lemma \ref{lem:2clique} shows that the above projection idea can be
used for other cliques and dense subgraphs.

The optimality and duality theory for the \SDP relaxations
is presented in Section \ref{sect:dualnonlin}.
We show that strict feasibility holds for the dual if the underlying
graph for the primal problem is connected.

Our primal-dual interior/exterior-point (p-d i-p) algorithm 
is derived in Section
\ref{sect:pdipNonlin}. We include a heuristic for obtaining a
strictly feasible starting point.
The algorithm uses a crossover technique, i.e. we use the affine
scaling step without backtracking once we get a sufficiently
large decrease in the duality gap.

We then continue with the numerical tests in Section \ref{sect:numeric}.
Concluding remarks are given in Section \ref{sect:concl}.

\section{\SNLt Problem Formulation}
\label{sect:formulation}
Let the $n$ unknown (sensor) points be $p^1, p^2, \ldots, p^n \in
\R^r$, $r$ the embedding dimension;
and let the $m$
known (anchor) points be $a^1, a^2, \ldots ,a^m \in \R^r$. Let
$X^T = [p^1, p^2, \ldots, p^n]$, and $A^T = [a^1, a^2, \ldots ,a^m]$.
We identify $a^i$ with $p^{n+i}$, for $i=1, \ldots, m$, 
and sometimes treat these as unknowns. We now define
\beq
\label{eq:Pdefn}
P^T:= 
\left( p^1, p^2, \ldots, p^n, a^1, a^2, \ldots ,a^m \right) = 
\left( p^1, p^2, \ldots, p^n, p^{n+1}, p^{n+2}, \ldots ,p^{n+m} \right) = 
\left( X^T \, A^T \right).
\eeq
Note that we can always translate all the sensors 
and anchors so that the anchors are
centered at the origin, i.e. $A^T \leftarrow A^T -\frac 1mA^Tee^T$ yields
$A^Te=0$. We can then translate them all back at the end.
In addition, we assume that there are a sufficient number of anchors 
so that the problem cannot be realized in a smaller embedding dimension.
Therefore, to avoid some special trivial cases, we assume the following.
\begin{assumption}
\label{assumpt:centered}
The number of sensors and anchors, and the embedding dimension satisfy
\[
n > > m > r, ~ A^Te=0, \mbox{ and }
A \mbox{ is full column rank}.
\]
\end{assumption}

Now define $\left(\NN_e, \NN_u, \NN_l\right)$, respectively,
to be the index sets of specified
(distance values, upper bounds, lower bounds), respectively,
of the distances $d_{ij}$ between pairs of nodes from $\{p^i\}_1^n$
(sensors);
and let $\left(\MM_e, \MM_u, \MM_l\right)$, denote the same for
distances between a node from $\{p^i\}_1^n$ (sensor) and a node
from $\{a^k\}_1^m$ (anchor).
Define (the partial Euclidean Distance Matrix) $E$ with elements
\[
E_{ij}=  \left\{ \begin{array}{ll}
  d_{ij}^2  & \mbox{if    }~~ ij \in \NN_e \cup \MM_e  \\
  \|p^{i}-p^{j}\|^2= \|a^{i-n}-a^{j-n}\|^2  & \mbox{if    }~~ i,j > n\\
  0       &   \mbox{otherwise}.
\end{array} \right.
\]
The underlying graph is
\beq
\label{eq:graphVE}
\GG=(\VV,\EE),
\eeq
with node set $\VV =\{1,\ldots,m+n\}$ and edge set 
$\EE =\NN_e \cup \MM_e\cup \{ij : i,j > n\}$.
Note that the subgraph induced by the anchors (the nodes with $j>n$) is
complete, i.e. the set of anchors forms a clique in the graph.
Similarly, we define the matrix of (squared) upper distance bounds $U^b$ and
the matrix of (squared) lower distance bounds $L^b$ for $ij \in \NN_u \cup
\MM_u$ and $\NN_l \cup \MM_l$, respectively.

Our first formulation for finding the sensor locations $p^j, j\leq n,$ is the
feasibility question for the constraints:
\beq
\label{eq:prelEDM1}
(\SNL_F)
  \begin{array}{rlcc}
        & \|p^i-p^j\|^2 = E_{ij} & \forall (i,j) \in \mathcal{N}_e
& \left(n_e=\frac {|\mathcal{N}_e|}{2} \right)\\
                          & \|p^i-a^k\|^2 = E_{ik} &  \forall (i,k) \in
\mathcal{M}_e & \left(m_e=\frac {|\mathcal{M}_e|}{2} \right)\\
        & \|p^i-p^j\|^2 = E_{ij} &  \forall i,j > n 
            & \left(\mbox{\underline{anchor-anchor}}\right)\\
        & \|p^i-p^j\|^2 \leq U^b_{ij} &  \forall (i,j) \in \mathcal{N}_u
& \left(n_u=\frac {|\mathcal{N}_u|}{2} \right)\\
                          & \|p^i-a^k\|^2 \leq U^b_{ik} &  \forall (i,k)
\in \mathcal{M}_u & \left(m_u=\frac {|\mathcal{M}_u|}{2} \right)\\
                          & \|p^i-p^j\|^2 \geq L^b_{ij} &  \forall (i,j)
\in \mathcal{N}_l & \left(n_l=\frac {|\mathcal{N}_l|}{2} \right)\\
                          & \|p^i-a^k\|^2 \geq L^b_{ik} &  \forall (i,k)
\in \mathcal{M}_l & \left(m_l=\frac {|\mathcal{M}_l|}{2}\right) \\
  \end{array}
\eeq
Note that the first three and the last two sets of constraints are
quadratic, nonconvex, constraints.
We added the anchor-anchor distances to emphasize that these are not
special and can be treated in the same way as the other distances.

The above can also be considered as a
{\em Graph Realization Problem}, i.e. we are given an
incomplete, undirected, edge-weighted simple graph. The
node set is the set of of sensors and anchors. The weights on the edges
are the squared distance between two nodes,
not all known and possibly inaccurate.
A Realization of $G$ in $\Re^r$ is a mapping of nodes 
into points $p^i$ in $\Re^r$ with squared distances given by the
weights.

Let $W_p, W_{pa},W_{a}$ be
weight matrices for the sensor-sensor, sensor-anchor, anchor-anchor,
distances respectively.
For example, they simply could be $0,1$ matrices to indicate
when an exact distance is unknown or known. 
Or a weight could be used
to verify the confidence in the value of the distance.
The weights in $W_a$ correspond to anchor-anchor distances and are 
{\em large}, since these distances are known.
In the literature,
these anchor-anchor distances are considered as constants in the
problem. We emphasize that they are equivalent to
the other distances, and that
the \SNL problem is a special case of the \EDM problem.
If there is noise in
the data, the exact model \eqref{eq:prelEDM1} can be infeasible. Therefore,
we can minimize the weighted least squares error.
\beq
\label{eq:edm}
(\SNL_{LS}) \qquad
  \begin{array}{rl}
      \mbox{min } f_1(P) := & \frac 12\displaystyle\sum_{(i,j) \in
\mathcal{N}_e}(W_p)_{ij}(\|p^i-p^j\|^2 - E_{ij})^2 \\
                          & + \frac 12\displaystyle\sum_{(i,k) \in
\mathcal{M}_e}(W_{pa})_{ik}(\|p^i-a^k\|^2 - E_{ik})^2 \\
                          & \left(+ \frac 12\displaystyle\sum_{i,j > n 
             }(W_a)_{ij}(\|p^i-p^j\|^2 - E_{ij})^2 \right)\\
       \mbox{subject to } & \|p^i-p^j\|^2 \leq U^b_{ij} \quad \forall (i,j)
\in \mathcal{N}_u \quad\left(n_u=\frac {|\mathcal{N}_u|}{2} \right)\\
                          & \|p^i-a^k\|^2 \leq U^b_{ik} \quad \forall (i,k)
\in \mathcal{M}_u \quad\left(m_u=\frac {|\mathcal{M}_u|}{2} \right)\\
                          & \|p^i-p^j\|^2 \geq L^b_{ij} \quad \forall (i,j)
\in \mathcal{N}_l \quad\left(n_l=\frac {|\mathcal{N}_l|}{2} \right)\\
                          & \|p^i-a^k\|^2 \geq L^b_{ik} \quad \forall (i,k)
\in \mathcal{M}_l \quad\left(m_l=\frac {|\mathcal{M}_l|}{2}\right)\\
        &\left( \|p^i-p^j\|^2 = E_{ij} \quad  \forall i,j > n\right).
  \end{array}
\eeq
This is a {\em hard} problem to solve due to the nonconvex
objective and constraints. We again included the anchor-anchor
distances within brackets both in the objective and constraints. This is
to emphasize that we could treat them with large weights in the
objective or as holding exactly without error in the constraints.

\section{Distance Geometry}
\label{sect:distgeom}
The geometry for \EDM has been studied in e.g.
\cite{sch35,gow85,hwlt91,MR14:889k}, and more recently, in e.g.
\cite{AlKaWo:97},\cite{homwolkA:04}.
Further theoretical properties can be found in e.g.
\cite{MR96a:15025,MR88k:15023,gow85,hwlt91,jt95,infieldsLaur:97,sch35,infieldsLaur:97}.
Since we emphasize that the \EDM theory can be used to solve the \SNL, we
now include an overview of the tools needed for \EDM. In particular, we
show the relationships between \EDM and \SDP.


\subsection{Linear Transformations and Adjoints Related to \EDMt}
(We use the notation from \cite{KrPiWo:06}. We include it here for
completeness.)
We work in spaces of real  matrices, $\MM^{s\times t}$, equipped with
the trace inner-product $\ip{A}{B}=\trace A^TB$ and induced
Frobenius norm $\|A\|_F^2=\trace A^TA$. For a given $B \in \Sn$, the
space of $n \times n$ real symmetric matrices, the linear transformation
$\diag(B) \in \R^n$ denotes the diagonal of $B$; for $v\in \R^n$, the adjoint
linear transformation is the diagonal matrix $\diag^*(v)=\Diag(v) \in \Sn$.
We now define several linear operators on $\Sn$.
(A collection of linear transformations, adjoints and properties
are given in the appendices.)
\beq
\label{eq:DKdef}
\begin{array}{rclrcl}
\DD_e(B) & := &  \diag (B)\,e^T + e \, \diag (B)^T, \qquad
\KK(B) & := &  \DD_e (B) - 2B,
\end{array}
\eeq
where $e$ is the vector of ones.  The adjoint linear operators are
\beq
\label{eq:DKadj}
\DD^*_e(D)= 2 \Diag(De), \qquad
\KK^*(D)=2(\Diag(De)-D).
\eeq
By abuse of notation we allow $\DD_e$ to act on $\R^n$:
\[
\DD_e(v) =   ve^T + ev^T,~ v\in \R^n.
\]
The linear operator $\KK$ maps the cone of
positive semidefinite matrices (denoted \SDP) onto the
cone of Euclidean distance matrices (denoted \EDM\@), i.e.
$\KK(\SDP)=\EDM$.  This allows us to
change problem \EDMC into a \SDP problem.
%

We define the linear transformation $\sblk_{i}(S)=S_{i} \in
{\mathcal S}^t$, on $S\in \Sn$,
that pulls out the $i$-th diagonal block of the matrix $S$
of dimension $t$. (The values of $t$ and $n$ can change and will
be clear from the context.) The adjoint $\sblk^*_{i}(T)=\sBlk_i(T)$,
where $T\in {\mathcal S}^t$, constructs a symmetric matrix of
suitable dimensions with all elements zero expect for the $i$-th
diagonal block given by $T$.

Similarly, we define  the linear transformation
$\sblk_{ij}(G)=G_{ij}$, on $G\in \Sn$, that pulls out the $ij$
block of the matrix $G$ of dimension $k\times l$ and multiplies it
by $\sqrt 2$. (The values of $k$, $l$, and $n$ can change and will
be clear from the context.) The adjoint
$\sblk^*_{ij}(J)=\sBlk_{ij}(J)$, where
$J\in \MM^{k\times l} \cong  \R^{kl}$,
constructs a symmetric matrix that has all elements zero expect
for the block $ij$ that is given by $J$ multiplied by $\frac
1{\sqrt 2}$, and for the block $ji$ that is given by
$J^T$multiplied by $\frac 1{\sqrt 2}$.
The multiplication by
$\sqrt 2$ (or $\frac 1{\sqrt 2}$) guarantees that the mapping is
an isometry.
We consider $J \in \MM^{k\times l}$ to be a $k \times l$ matrix and
equivalently $J \in \R^{kl}$ is a vector of length $kl$ with the
positions known.

\subsection{Properties of Transformations}
\begin{lem}
\label{lem:nullrange}
(\cite{homwolkA:04})
Define the linear operator on $\Sn$ by
\[
\offDiag (S) = S - \Diag(\diag(S)).
\]
Let $J:=I-\frac 1n ee^T$. Then, the following holds.
\begin{itemize}
\item[$\bullet$] The nullspace $\NN(\KK)$ equals the range $\RR(\DD_e)$.
\item[$\bullet$] The range $\RR(\KK)$ equals the hollow subspace of
$\Sn$, denoted $S_H:=\{D\in\Sn: \diag(D)=0\}$.

\item[$\bullet$] The range $\RR(\KK^*)$ equals the centered subspace of
$\Sn$, denoted $S_c:=\{B\in \Sn: Be=0\}$.
\item[$\bullet$] The Moore-Penrose generalized inverse 
$\KK^\dagger (D) = -\frac 12 J\left(\offDiag(D)\right)J$.
\epr
\end{itemize}
\end{lem}
\begin{cor}(\cite{KrPiWo:06})
\begin{enumerate}
\item
\label{item:decomp1}
Let $S_D$ denote the subspace of diagonal matrices in $\Sn$. Then
\[
\begin{array}{lcl}
S_c = \NN(D_e^*) = \RR(\KK^*) = \RR(\KK^\dagger)& \perp & \NN(K)=\RR(\DD_e)  \\
S_H = \RR(K) = \NN(\DD_e) & \perp & S_D=\NN(K^*)=\RR(\DD_e^*).
\end{array}
\]
\item
\label{item:decomp2}
Let $\pmatrix{V & \frac 1{\sqrt{n}} e}$ be an $n \times n$ orthogonal
matrix.
Then
\[
Y \succeq 0 \iff
Y=V \hat Y V^T + \DD_e(v) \succeq 0, \mbox{ for some } \hat Y \in
\Snm, v \in \R^n.
\]
\end{enumerate}
\epr
\end{cor}

Let $B=PP^T$. Then
\beq
\label{eq:d2B}
D_{ij}=\|p^i-p^j\|^2 = (\diag (B) e^T+ e \diag(B)^T -2B)_{ij}
      =\left(\KK(B)\right)_{ij},
\eeq
i.e. the \EDM $D=(D_{ij})$ and the points $p_i$ in $P$
are related by $D= \KK(B)$, see \eqref{eq:DKdef}.
\begin{lem}(\cite{KrPiWo:06})
Suppose that $0 \preceq B \in \Sn$. Then $D=\KK(B)$ is EDM.
\epr
\end{lem}

\section{\SDPt Relaxations of \SNLt based on \EDMt Model}
\label{sect:relax}

We first study the \SDP relaxation used in the recent series of papers
on \SNL, e.g.
\cite{BiswasYe:04,biswasliangtohwangye,SoYe:05,biswasliangtohye:05,Jin:05}.
(See \eqref{eq:currentsdp} and Section \ref{sect:currentproj} below.)
This relaxation starts by treating the anchors
distinct from the sensors.
We use a different derivation and model the problem based on 
classical \EDM theory, and show its equivalence with the current 
\SDP relaxation. By viewing the \SNL problem as an \EDM problem, 
we obtain several interesting results, e.g. clique reduction, and 
a geometric interpretation on how to estimate sensor positions from the
\SDP relaxation optimum.

\subsection{Connections from Current \SDPt Relaxation to \EDMt}
\label{sect:relaxcurrent}
Let $Y= XX^T$.  Then the current
\SDP relaxation for the feasibility problem for \SNL uses 
\beq
\label{eq:relaxZXY}
Y \succeq  XX^T, \mbox{ or equivalently, }
Z_s=\pmatrix{I_r & X^T \cr X & Y } \succeq 0.
\eeq
 This is in combination with the constraints 
\beq
\label{eq:currentsdp}
\begin{array}{rcl}
\trace \pmatrix{ 0 \cr e_i -e_j }
\pmatrix{ 0 \cr e_i -e_j }^T Z_s &=&
 E_{ij}, \quad \forall ij\in \NN_e  \\
\trace \pmatrix{ -a_k \\ e_i  }
 \pmatrix{-a_k \cr e_i  }^T Z_s &=&
 E_{ij}, \quad \forall ij \in \MM_e, i<j=n+k.
\end{array}
\eeq

\subsubsection{Reformulation using Matrices}
We use the matrix lifting or linearization
$\bar Y:= PP^T = \pmatrix{XX^T & XA^T \cr AX^T &  AA^T}$
and $Z:=[I;P][I;P]^T=\pmatrix{I & P^T\cr P & \bar Y}$.
The dimensions are:
\[
X \in \MM^{n,r}; \quad A \in \MM^{m,r}; \quad P \in \MM^{m+n,r};
\quad \bar Y \in \Ss^{m+n}; \quad Z \in \Ss^{m+n+r}.
\]
Adding the hard quadratic constraint $\bar Y= PP^T$ allows us to replace the
quartic objective function in \SNLLS with a quadratic function.
We can now reformulate \SNL using matrix notation to get the
equivalent \EDM problem
\beq
\label{eq:edm2}
(\SNL_M) \qquad
  \begin{array}{crcl}
      \min &   f_2(\bar Y) :=  \frac{1}{2}{\| W \circ (\KK(\bar Y) - E)
\|}_F^2 \\
\mbox{subject to } & g_u(\bar Y):=H_u\circ(\KK(\bar Y) - \bar U^b) &\leq& 0 \\
                          & g_l(\bar Y) :=H_l \circ (\KK(\bar Y) - \bar L^b) &
\geq& 0\\
                          & \bar Y-PP^T&=&0\\
  & \left(\sblk_2 (\KK(\bar Y))\right. &=& \left. \KK(AA^T) \right),
  \end{array}
\eeq
where $W\in \Ss^{n+m}$ is the weight matrix having a
positive ij-element if $(i,j)\in \mathcal{N}_e \cup
\mathcal{M}_e \cup \{(ij): i,j>n\}$, $0$ otherwise. 
$H_u, H_l$ are $0,1$-matrices where
the ij-th element equals $1$ if an upper (resp. lower) bound exists; and
it is  0 otherwise. By abuse of notation, 
we consider the functions $g_u,g_l$ as
implicitly acting on only the nonzero components in the upper triangular
parts of the matrices that result from the Hadamard
products with $H_u,H_l$, respectively.
We include in brackets the constraints on the clique corresponding to
the anchor-anchor distances.
\begin{rem}
The function $f_2(\bar Y)=f_2(PP^T)$, and it is clear
that $f_2(PP^T)=f_1(P)$ in \eqref{eq:edm}.
Note that the functions $f_2,g_u,g_l$ act only on
$\bar Y$ and the locations of the anchors and sensors is completely hiding
in the {\em hard}, nonconvex quadratic constraint
$\bar Y= PP^T = \pmatrix{XX^T & XA^T \cr AX^T &  AA^T}$.
The problem \SNLM is a linear least squares problem with nonlinear
constraints. The objective function is generally
underdetermined.  
This can result in ill-conditioning problems, e.g. \cite{MR2059740}.
Therefore, reducing the number of variables helps with stability.
\end{rem}

\subsubsection{\SDPt Relaxation of the Hard Quadratic Constraint}
\label{sect:relaxhard}
We now consider the hard quadratic constraint in \eqref{eq:edm2} 
\beq
\label{eq:YPP}
\bar Y= 
 \pmatrix{\bar Y_{11} & \bar Y_{21}^T \cr \bar Y_{21} &  AA^T}=
PP^T = \pmatrix{XX^T & XA^T \cr AX^T &  AA^T},
\eeq
where $P$ is defined in \eqref{eq:Pdefn}.  We study the standard current
semidefinite relaxation in \eqref{eq:currentsdp} with
\eqref{eq:relaxZXY}, or equivalently with $\bar Y \succeq PP^T$. 
We show that this is equivalent to the simpler $\bar Y \succeq 0$.
We include details on problems and
weaknesses with the relaxation. We first present
several lemmas. We start with the following well known
result. We include a proof for completeness.
\begin{lem}
\label{eq:lemrangeY}
Suppose that the partitioned symmetric
matrix $\pmatrix{Y_{11} & Y_{21}^T \cr Y_{21} & AA^T } \succeq 0$.
Then $Y_{21}^T=XA^T$, with $X= Y_{21}^TA(A^TA)^{-1}$.
\end{lem}
\bpr
Let $A=U\Sigma_r V^T$ be the compact singular value decomposition,
$0 \prec \Sigma_r \in \Ss^r$. And, suppose that $\pmatrix{ U & \bar U}$ is
an orthogonal matrix. Therefore, the range spaces
$\RR(U)=\RR(A)$ and the nullspace $\NN(A^T)=\RR(\bar U)$.
Consider the nonsingular congruence
\[
\begin{array}{rcl}
0 
&\preceq &
\pmatrix{ I & 0 \cr 0 & \pmatrix{ U &  \bar U}}^T
\pmatrix{Z & Y_{21}^T \cr Y_{21} & AA^T } 
\pmatrix{ I & 0 \cr 0 & \pmatrix{ U &  \bar U}}
\\&= &
\pmatrix{Z & Y_{21}^T\pmatrix{ U &  \bar U} \cr \pmatrix{ U &  \bar U}^T Y_{21} & 
                       \pmatrix{ \Sigma_r^2 &  0 \cr 0 & 0}}. 
\end{array}
\]
This implies that $Y_{21}^T\bar U = 0$. This in turn means that
$\NN(Y_{21}^T) \supset \NN(A^T)$, or equivalently, $\RR(A) \supset
\RR(Y_{21})$. Note that the orthogonal projection onto $\RR(A)$ is
$A (A^TA)^{-1}A^T$. Therefore, 
$Y_{21}^T= Y_{21}^TA(A^TA)^{-1}A^T =\left(Y_{21}^TA(A^TA)^{-1}\right)A^T$,
i.e. we can choose $X= Y_{21}^TA(A^TA)^{-1}$.
\epr

In the recent literature, e.g.
\cite{MR2191577,BiswasYe:04,Jin:05},
it is common practice to relax the hard constraint
\eqref{eq:YPP} to a tractable semidefinite
constraint, $\bar Y \succeq PP^T$, or equivalently, $\bar Y_{11} \succeq XX^T$
with $\bar Y_{21}=AX^T$. The following lemma presents several
characterizations for the resulting feasible set.
\begin{lem}
\label{lem:equivcontr}
Let  $A=U\Sigma_r V^T$ be the compact singular value decomposition of
$A$, and let
$P,\bar Y$ be partitioned as in \eqref{eq:Pdefn},\eqref{eq:YPP}, 
\[
P = \pmatrix{ P_1 \cr P_2}, \quad
\bar Y= \pmatrix{\bar Y_{11} & \bar Y_{21}^T \cr \bar Y_{21} &  \bar Y_{22}}.
\]
Define the semidefinite relaxation of the hard quadratic constraint
\eqref{eq:YPP} as:
\beq
\label{eq:bigG}
G(P,\bar Y):= PP^T-\bar Y \preceq 0, \quad \bar Y_{22}=AA^T,
  \quad P_2=A.
\eeq
By abuse of notation, we allow $G$ to act on spaces of different
dimensions. Then we get
the following equivalent representations of the
corresponding feasible set $\FF_G$.
\begin{subequations}
\beq
\label{eq:equivconstr0}
 \FF_G=\left\{(P,\bar Y):  G(P,\bar Y) \preceq 0, \bar Y_{22}=AA^T,
       P_2=A \right\}
\eeq
\beq
\label{eq:equivconstr1}
 \FF_G=\left\{(P,\bar Y):  G(X,Y) \preceq 0, 
  \bar Y_{11}=Y, \bar Y_{21}=AX^T,\bar Y_{22}=AA^T, P=\pmatrix{X\cr A} \right\}
\eeq
\beq
\label{eq:equivconstr2}
\begin{array}{cc}
 \FF_G=\left\{(P,\bar Y):  
Z= \pmatrix{Z_{11} & Z_{21}^T \cr Z_{21}  & Z_{22} }
\succeq 0, 
   \bar Y= \pmatrix{ I_n & 0 \cr 0 & A } Z \pmatrix{ I_n & 0 \cr 0 & A }^T,
       \right. \\
       \left. \qquad\qquad\qquad    Z_{22}=I,
           X= Z_{21}^T, P=\pmatrix{X\cr A}
                \right\}
\end{array}
\eeq
\beq
\label{eq:equivconstr3}
 \FF_G=\left\{(P,\bar Y):  \bar Y \succeq 0, 
      \bar Y_{22}=AA^T, X= \bar Y_{21}^TA(A^TA)^{-1},
               P=\pmatrix{X\cr A} \right\}
\eeq
\beq
\label{eq:equivconstr4}
\begin{array}{cc}
 \FF_G=\left\{(P,\bar Y):  
Z= \pmatrix{Z_{11} & Z_{21}^T \cr Z_{21}  & Z_{22} }
\succeq 0, 
   \bar Y= \pmatrix{ I_n & 0 \cr 0 & U } Z \pmatrix{ I_n & 0 \cr 0 & U }^T,
\right. \\ 
   \left. \qquad \qquad\qquad    Z_{22}=\Sigma_r^2,
           X= \bar Y_{21}^TA(A^TA)^{-1}=Z_{21}^T\Sigma_r^{-1}V^T, P=\pmatrix{X\cr A}
                \right\}.
\end{array}
\eeq
\end{subequations}
Moreover, the function $G$ is convex in the L\"{o}wner (semidefinite) partial order;
and the feasible set $\FF_G$ is a closed convex set.
\end{lem}
\bpr
Recall that the cone of positive semidefinite matrices is self-polar.
Let $Q\succeq 0$ and $\phi_Q(P)=\trace QPP^T$.
Convexity of $G$ follows from positive semidefiniteness of the Hessian
$\nabla^2 \phi_Q(P)=I \otimes Q$, where $\otimes$ denotes the Kronecker
product.

In addition,
\[
0 \succeq G(P,\bar Y)= PP^T-\bar Y =
\pmatrix{
   XX^T - \bar Y_{11} & XA^T - \bar Y_{21}^T  \cr
    AX^T - \bar Y_{21} & 0
}
\]
holds if and only if
\[
0 \succeq G(X,\bar Y_{11})= XX^T-\bar Y_{11}, \mbox{  and  }
    AX^T - \bar Y_{21} = 0.
\]
This shows the equivalence with
\eqref{eq:equivconstr1}.
A Schur complement argument, with $\bar Y_{11}=Y$, shows the equivalence with
$\pmatrix{Y & X \cr X^T & I_r } \succeq 0$, i.e. with the set in
\eqref{eq:equivconstr2}.
The equivalence with
\eqref{eq:equivconstr3}
follows from Lemma \ref{eq:lemrangeY}.

To show the equivalence with the final expression
\eqref{eq:equivconstr4}, we note
that $\bar Y \succeq 0, \bar Y_{22}=AA^T$, implies that there is no
strictly feasible $\bar Y \succ 0$. Therefore, we project the feasible
set onto the {\em minimal cone} or face (see \cite{bw1}). This yields the 
minimal face that contains the feasible set of $\bar Y$, i.e.
\beq
\label{eq:minimalface}
\bar Y= \pmatrix{ I_n & 0 \cr 0 & U } Z_w \pmatrix{ I_n & 0 \cr 0 & U }^T,
Z_w \succeq 0, Z_w\in \Ss^{n+r}.
\eeq
The result follows since the constraint $\bar Y_{22}=AA^T$ holds if and only if
$Z_w$ is blocked as
$Z_w:= \pmatrix{Y& W^T \cr W & \Sigma_r^2 } \succeq 0$. 
(More simply, one can show the equivalence of \eqref{eq:equivconstr4} with
\eqref{eq:equivconstr2} by using the compact singular value decomposition of
$A$. However, the longer proof given above emphasizes that the reduction comes
from using a projection to obtain the Slater constraint qualification.)
\epr
The above Lemma \ref{lem:equivcontr}
shows that we can treat the set of anchors as a set
of sensors for which all the distances are known, i.e. the set
of corresponding nodes is a clique. The fact that we have a
clique and the diagonal $m \times m$ block $AA^T$ in $\bar Y$ is rank deficient, $r<m$,
means that the Slater constraint qualification, $\bar Y \succ 0$, cannot hold.
Therefore, we can project onto the minimal cone containing the feasible
set and thus reduce the size of the problem, see Lemma
\ref{lem:equivcontr}, \eqref{eq:equivconstr4}, i.e. the variable
$\bar Y \in \Ss^{n+m}$ is reduced in size to $Z \in \Ss^{n+r}$.
The reduction can be done by using any point in the relative interior of
the minimal cone, e.g. any feasible point of maximum rank. The equivalent
representations in \eqref{eq:equivconstr2} and in \eqref{eq:equivconstr4}
illustrate this.

\subsubsection{Current \SDPt Relaxation using Projection onto Minimal Cone}
\label{sect:currentproj}
The above reduction to $Y$ in Lemma \ref{lem:equivcontr},
\eqref{eq:equivconstr1},
allows us to use the smaller dimensional semidefinite constrained
variable
\beq
\label{eq:Zs}
Z_s=
\pmatrix{I & X^T \cr X & Y}
\succeq 0 \in \Ss^{n+m}, \quad
\bar Y_{11}=Y, \bar Y_{21}=AX^T.
\eeq
This is what is introduced in e.g. \cite{BiswasYe:04}.

\begin{rem}
\label{rem:onto}
Note that the mapping $Z_s=Z_s(X,Y): \MM^{n,r} \times \Ss^n \rightarrow
\Ss^{n+r}$ is {\em not} onto. This means that the Jacobian of the
optimality conditions cannot be full rank, i.e. this formulation
introduces instability into the model. A minor modification corrects
this, i.e.  the $I$ constraint is added explicitly.
\[
Z=
\pmatrix{Z_{11} & Z_{21}^T \cr Z_{21} & Z_{22}}
\succeq 0, \quad
Z_{11}=I, \bar Y_{11}=Z_{22},  \bar Y_{21}=AZ_{21}^T.
\]
\end{rem}

To develop the model for computations, we introduce the following
notation.
\[
x:= \kvec\left( \sblk_{21} \pmatrix{ 0 & X^T \cr X & 0}\right)=
\sqrt{2} \kvec(X), \quad  y:=\svec(Y),
\]
where we add $\sqrt{2}$ to the definition of $x$ since $X$ appears together
with $X^T$
in $Z_s$ and implicitly in $\bar Y$, with $Y_{21}=AX^T$.
We define the following matrices and linear transformations:
\[
\begin{array}{rcl}
    {\ZZ}^x_s(x):=\sBlk_{21}(\Mat(x)), &
                    {\ZZ}^y_s(y):=\sBlk_2(\sMat(y)),\\
       {\ZZ}_s(x,y):=\ZZ^x_s(x)+\ZZ^y_s(y),
                &
    Z_s  := \sBlk_1(I)+{\ZZ}_s(x,y), \\
\end{array}
\]
\[
\begin{array}{rcl}
\YY^x(x):=\sBlk_{21}(A\Mat(x)^T), &
\YY^y(y):=\sBlk_1(\sMat(y))\\
\YY(x,y):= \YY^x(x)+ \YY^y(y), & \bar Y:=\sBlk_2(AA^T)+\YY(x,y).
\end{array}
\]
\[
\begin{array}{rcl}
\vspace{.02in}
\bar E:= W \circ \left[ E- \KK(\sBlk_2(AA^T))\right],\\
\vspace{.02in}
\bar U^b:=  H_u \circ \left[\KK(\sBlk_2(AA^T)) - U^b\right],\\
\bar L^b:= H_l\circ \left[L^b-\KK(\sBlk_2(AA^T))\right].
\end{array}
\]
By abuse of notation, we let the functions ${\ZZ}^x_s, \ldots, \YY$ act
directly on the matrices $X,Y$. The meaning will be clear from the
context.

The unknown matrix $\bar Y$ in \eqref{eq:edm2} is equal to
$\YY(x,y)$ with the additional constant in the
$2,2$ block, i.e. our unknowns are the vectors $x,y$ which are used to
build $\bar Y$ and $Z_s$.
Using this notation we can introduce the following vector
form of the relaxation of \eqref{eq:edm2}.
\beq \label{eq:relax1}
(SNL_{MV}) \qquad
  \begin{array}{crcl}
      \min  &  f_3(x,y) :=  \frac{1}{2}{\| W \circ (\KK(\YY(x,y))) - \bar E
\|}_F^2 \\
  \mbox{subject to } & g_u(x,y):=H_u \circ \KK(\YY(x,y)) - \bar U^b &\leq& 0
\\
                   & g_l(x,y):=\bar L^b -H_l \circ \KK(\YY(x,y))& \leq& 0 \\
                & \sBlk_1(I)+{\ZZ}_s(x,y)&\succeq& 0.
  \end{array}
\eeq
As above, we consider the functions $g_u,g_l$ as implicitly
acting only on the
nonzero parts of the upper triangular part of the
matrix that results from the Hadamard
products with $H_u,H_l$, respectively.

\subsubsection{\SDPt Formulation Using \EDMt}

The equivalent representations of the feasible set given in
Lemma \ref{lem:equivcontr}, in particular by \eqref{eq:equivconstr4}, 
show that \SNL is an \EDM problem $D=\KK(\bar Y)$, with the additional upper
and lower bound constraints as well as the block constraint
$\sBlk_2(D)=\KK(AA^T)$, or equivalently, $\sBlk_2(\bar Y)=AA^T$. 

\begin{rem}
Suppose that we can increase the size of the clique
containing the anchor nodes by adding sensor nodes where the distances are
exactly known. Then these sensor nodes can be treated as anchor nodes, though
their position is unknown.
\end{rem}

We can now obtain an equivalent relaxation for \SNL by using the \EDM
\eqref{eq:edm2} and replacing the hard quadratic constraint with the simpler
semidefinite constraint $\bar Y \succeq 0$.
We then observe that the Slater constraint qualification (strict feasibility)
fails. Therefore, we can project onto the minimal cone, i.e. onto the minimal
face of the \SDP cone that contains the feasible set. see \cite{bw1,AlKaWo:97}.
Let 
\beq
\label{eq:Us}
U_A=  \pmatrix{ I_n & 0 \cr 0 & A }, \quad
U_s=  \pmatrix{ I_n & 0 \cr 0 & U }, 
\mbox{ where } A=U\Sigma_r V^T.
\eeq
Recall that 
  $\sblk_2 (\KK(U_s Z U_s^T)) =  \KK(AA^T)$ is equivalent to
$Z_{22}=\Sigma^2_r$.
We get two \SDP relaxation that are equivalent to \SNLMV:
\beq
\label{eq:snl2s}
(\SNL_{\scriptsize{\EDM s}}) 
  \begin{array}{crcl}
      \min &   f_3(Z) :=  \frac{1}{2}{\| W \circ (\KK(U_s Z U_s^T) - E)
\|}_F^2 \\
\mbox{subject to } & H_u \circ (\KK(U_s Z U_s^T) - U^b) \leq 0 \\
                          & H_l \circ (\KK(U_s Z U_s^T) - L^b) \geq 0\\
  & Z_{22}=\Sigma^2_r  \\
  &   Z\succeq 0.
  \end{array}
\eeq
and
\beq
\label{eq:snl2A}
(\SNL_{\scriptsize{\EDM A}}) 
  \begin{array}{crcl}
      \min &   f_3(Z) :=  \frac{1}{2}{\| W \circ (\KK(U_A Z U_A^T) - E)
\|}_F^2 \\
\mbox{subject to } & H_u \circ (\KK(U_A Z U_A^T) - U^b) \leq 0 \\
                          & H_l \circ (\KK(U_A Z U_A^T) - L^b) \geq 0\\
  & Z_{22}=I_r  \\
  &   Z\succeq 0.
  \end{array}
\eeq
\begin{rem}
Note that we do not substitute the constraint on $Z_{22}$ into $Z$, but leave
it explicit. Though this does not change the feasible set, it does change the
stability and the dual. This can be compared to the \SDP relaxation for the
Max-Cut problem with constraint that the diagonal of $X$ is all ones,
$\diag X = e$ and $X\succeq 0$. However, one does not substitute for the
diagonal and rewrite the semidefinite constraint.
\end{rem}

\subsection{Clique Reductions using Minimal Cone Projection}

Now suppose that we have another clique of $p>r$ sensors 
where the exact distances are known and are used as constraints. 
Then there exists a matrix $\bar Y=PP^T$ that has a diagonal rank deficient 
$p \times p$ block.  Since all feasible points are found from elements
in the set $\bar Y + \NN(\KK)$, we conclude that for $p$ large
enough, the diagonal block remains rank deficient for all feasible $\bar
Y$, i.e. the Slater constraint qualification fails again, if the
corresponding distances are added as constraints.

We now see that we can again take advantage of the loss of the Slater
constraint qualification.
\begin{lem}
\label{lem:2clique}
Suppose that the hypotheses and definitions from
Lemma \ref{lem:equivcontr} hold; and
suppose that there exists a set of sensors, without loss of generality
$S_c:=\{p^{t+1}, \ldots, p^n\}$, so that the distances
$\|p^i-p^j\|$ are known for all $t+1\leq i,j\leq n$;
$a^k$,
i.e. the graph of the partial \EDM has two cliques, one clique
corresponding to the set of known anchors, and the other to the set of
sensors $S_c$.
Let $P,\bar Y$ be partitioned as
\[
P = \pmatrix{ P_1 \cr P_2\cr P_3}= \pmatrix{X \cr A}, \quad
\bar Y= \pmatrix{
          \bar Y_{11} & \bar Y_{21}^T  & \bar Y_{31}^T \cr
        \bar Y_{21} &  \bar Y_{22} & \bar Y_{32}^T \cr
        \bar Y_{31} &  \bar Y_{32} & \bar Y_{33}
          }=PP^T,
\]
where $P_i=A_i, i=2,3$, and $A_3=A$ corresponds to the known
anchors while $P_2=A_2$ corresponds to the clique of sensors and
$X=\pmatrix{P_1\cr P_2}$ corresponds to all the sensors.
Let the \EDM, $E=\KK(\bar Y)$, be correspondingly blocked
\[
E=\pmatrix{E_1 &\cdot  &\cdot  \cr \cdot &
E_2 &\cdot  \cr \cdot &\cdot &E_3},
\]
so that $E_3=\KK(AA^T)$ are the anchor-anchor squared distances, and
$E_2=\KK(P_2P_2^T)$ are
the squared distances between the sensors in the set $S_c$.
Let
\[
B=  \KK^\dagger (E_2).
\]
Then the following hold.
\begin{enumerate}
\item $Be=0$ and
\beq
\label{eq:YK}
 \bar Y_{22}= B +\bar y_2e^T +e \bar y_2^T  \succeq 0,
\mbox{ for some } \bar y_2 \in \RR(B)+\alpha e, \alpha \geq 0,
 \mbox{ with } \rank(\bar Y_{22}) \leq r.
\eeq
\item The feasible set $\FF_{G}$ in Lemma \ref{lem:equivcontr} can be
formulated as
\beq
\label{eq:subfeasset}
\begin{array}{cc}
 \FF_{G}:=\left\{(P,\bar Y):
Z= \pmatrix{Z_{11} & Z_{21}^T& Z_{31}^T \cr
        Z_{21}  & Z_{22} & Z_{32}^T  \cr
        Z_{31}  & Z_{32} & Z_{33} }
\succeq 0,
   \bar Y=
  \pmatrix{ I_t & 0& 0 \cr 0 & U_2 & 0\cr 0 &0 & U_s } Z
  \pmatrix{ I_t & 0& 0 \cr 0 & U_2 & 0\cr 0 &0 & U_s }^T,
\right. \\
   \left.  \qquad\qquad      Z_{33}=\Sigma_r^2,
           X= \pmatrix{Z_{31}\cr U_2Z_{32}}\Sigma_r^{-1}V_s^T,
      P=\pmatrix{X\cr A} \right\},
\end{array}
\eeq
or equivalently as
\beq
\label{eq:subfeassetA}
\begin{array}{cc}
 \FF_{G}=\left\{(P,\bar Y):
Z= \pmatrix{Z_{11} & Z_{21}^T& Z_{31}^T \cr
        Z_{21}  & Z_{22} & Z_{32}^T  \cr
        Z_{31}  & Z_{32} & Z_{33} }
\succeq 0,
   \bar Y=
  \pmatrix{ I_t & 0& 0 \cr 0 & U_2 & 0\cr 0 &0 & A } Z
  \pmatrix{ I_t & 0& 0 \cr 0 & U_2 & 0\cr 0 &0 & A }^T,
\right. \\
   \left.  \qquad\qquad      Z_{33}=I_r,
           X= \pmatrix{Z_{31}^T\cr U_2Z_{32}^T},
      P=\pmatrix{X\cr A} \right\},
\end{array}
\eeq
where $\hat B:=B+2ee^T= \pmatrix{U_2 & \bar U_2} \pmatrix{ D_2 & 0 \cr 0 &
0}
\pmatrix{U_2 & \bar U_2}^T$ is the orthogonal
diagonalization of $\hat B$, with $D_2 \in \Ss_+^{r_2}, r_2\leq r+1$.
\end{enumerate}
\end{lem}
\bpr
We proceed just as we did in Lemma \ref{lem:equivcontr}, i.e. we reduce the
problem by projecting onto a smaller
face in order to obtain the Slater constraint qualification.

The equation for $\bar Y_{22}$ for some $\bar y_2$, given in \eqref{eq:YK}, follows
from the nullspace characterization in Lemma \ref{lem:nullrange}.
Moreover, $\bar Y_{22}=P_2P_2^T$ implies that $\rank (Y_{22}) \leq r$, the
embedding dimension. And, $\bar Y_{22} \succeq 0, Be=0$ implies the inclusion
$\bar y_2 \in \RR(B)+\alpha e, \alpha \geq  0$.
Moreover, we can shift $\bar P_2^T=P_2^T-\frac 1{n-t}(P_2^Te)e^T$. Then
for $B=\bar P_2 \bar P_2^T$, we get $Be=0$, i.e. this satisfies
$B=\KK^\dagger(E_2)$ and $\rank(B) \leq r$. Therefore,
for any $Y=B+ye^T+ey^T
\succeq 0$, we must have $y = \alpha e, \alpha \geq 0$.
Therefore, $\hat B$ has the maximum rank, at most $r+1$, among all feasible matrices
of the form $0 \preceq Y \in B+\NN(\KK)$. $\hat B$ determines the
smallest face containing all such feasible $Y$.

Define the linear transformation
$H:\R^{n-t} \rightarrow \Ss^{n-t}$ by $H(y)=\bar Y_{22}+ye^T+ey^T$.
Let $\LL := \bar Y_{22} + \RR(\DD_e)$ and $\FF_e$ denote the smallest face of
$\Ss^{n-t}_+$ that contains $\LL$. Since $\hat B$ is a feasible
point of maximum rank, we get
\[
\hat B = B + \DD_e(\hat y_2) \in  \left(\LL \cap \relint \FF_e\right).
\]
Thus the face 
\[
\FF_e= 
\{ U_2ZU_2^T:  Z\in \Ss^{r_2}_+\}=
\{ Y \in \Ss^{n-t}_+:  \trace Y(\bar U_2 \bar U_2^T)=0 \}.
\]

Now, we expand
\[
\begin{array}{rcl}
\pmatrix{\bar Y_{11} & \bar Y_{21}^T  & \bar Y_{31}^T \cr
        \bar Y_{21} &  \bar Y_{22} & \bar Y_{32}^T \cr
        \bar Y_{31} &  \bar Y_{32} & \bar Y_{33}}
        &=&
  \pmatrix{ I_t & 0& 0 \cr 0 & U_2 & 0\cr 0 &0 & U_s }
     \pmatrix{Z_{11} & Z_{21}^T& Z_{31}^T \cr 
        Z_{21}  & Z_{22}   & Z_{32}^T  \cr 
        Z_{31}  & Z_{32} & \Sigma_r^2 }
  \pmatrix{ I_t & 0& 0 \cr 0 & U_2 & 0\cr 0 &0 & U_s }^T
        \\&=&
     \pmatrix{Z_{11} & Z_{21}^TU_2^T& Z_{31}^TU_s^T \cr 
        U_2Z_{21}  & U_2Z_{22}U_2^T   & U_2Z_{32}^TU_s^T  \cr 
       U_s Z_{31}  & U_sZ_{32}U_2^T & U_s\Sigma_r^2U_s^T }.
\end{array}\]
Therefore, $\bar Y_{21}^T=\pmatrix{Z_{31}U_s^T\cr U_2Z_{32}^TU_s^T}$.
Therefore, the expression for $Z_{33}$ and $X$ in \eqref{eq:subfeasset} 
follows from
equation \eqref{eq:equivconstr4} in Lemma \ref{lem:equivcontr}.
The result in \eqref{eq:subfeasset} can be obtained similarly or by
using the compact singular value decomposition of $A$.
\epr

\begin{rem}
The above Lemma \ref{lem:2clique} can be extended to sets of sensors 
that are not cliques, but have many known
edges. The key idea is to be able to use 
$(W_i\circ \KK)^\dagger (W_i\circ E_i)$ and to characterize the nullspace of
$W_i\circ \KK$. This is studied in a forthcoming paper.
\end{rem}
We can apply Lemma \ref{lem:2clique} to further reduce the \SDP relaxation. 
Suppose there are a group of sensors for which pairwise distances 
are all known. This should be a common occurrence,
since distances between sensors within radio range are all known. 
Without loss of generality, 
we assume the of sensors to be $\{p^{t+1}, \ldots, p^{n}\}$. 
Let $E_2$, $B=\KK^{\dag}(E_2)$, and $U_2$, be found using
Lemma \ref{lem:2clique} and denote
\beq
\label{eq:U2s}
U_{2A}:=  \pmatrix{ I_n & 0 &0 \cr 0 & U_2 &0\cr 0 & 0& A }, \quad 
U_{2s}:=  \pmatrix{ I_n & 0 &0 \cr 0 & U_2 &0\cr 0 & 0& U }. 
\eeq
In $\SNL_{\scriptsize{\EDM s}}$, 
we can replace $U_{s}$ with $U_{2s}$ and reach a reduced \SDP formulation. 
Similarly, for $\SNL_{\scriptsize{\EDM A}}$.
Furthermore, we may generalize to the $k$ clique cases 
for any positive integer $k$. We similarly define each $U_i, 
2\leq i \leq k$, and define
\beq
\label{eq:Uks}
U_{kA}=  \pmatrix{ I_n & 0 & \ldots &0&0 \cr 0 & U_2 & \ldots&0 & 0\cr \vdots &&\ddots&&\vdots \cr 0 & 0 &0 & U_k &0 \cr 0 &0&0&0& A}, \quad
U_{ks}=  \pmatrix{ I_n & 0 & \ldots &0&0 \cr 0 & U_2 & \ldots&0 & 0\cr \vdots &&\ddots&&\vdots \cr 0 & 0 &0 & U_k &0 \cr 0 &0&0&0& U}.
\eeq
Then we can formulate a reduced \SDP for $k$ cliques:
\beq
\label{eq:snl3A}
(\SNL_{\scriptsize{k-cliques-A}}) 
  \begin{array}{crcl}
      \min &   f_4(Z) :=  \frac{1}{2}{\| W \circ (\KK(U_{kA} Z U_{kA}^T) - E)
\|}_F^2 \\
\mbox{subject to } & H_u \circ (\KK(U_{kA} Z U_{kA}^T) - U^b) \leq 0 \\
                          & H_l \circ (\KK(U_{kA} Z U_{kA}^T) - L^b) \geq 0\\
  & Z_{kk}=I_r  \\
  &   Z\succeq 0
  \end{array}
\eeq
where $Z_{kk}$ is the last $r$ by $r$ diagonal block of $Z$.  
Similarly, we get
\beq
\label{eq:snl3s}
(\SNL_{\scriptsize{k-cliques-s}}) 
  \begin{array}{crcl}
      \min &   f_4(Z) :=  \frac{1}{2}{\| W \circ (\KK(U_{ks} Z U_{ks}^T) - E)
\|}_F^2 \\
\mbox{subject to } & H_u \circ (\KK(U_{ks} Z U_{ks}^T) - U^b) \leq 0 \\
                          & H_l \circ (\KK(U_{ks} Z U_{ks}^T) - L^b) \geq 0\\
  & Z_{kk}=I_r  \\
  &   Z\succeq 0
  \end{array}
\eeq

For a clique with $r_e$ sensors, a $U_i$ is constructed with $r_e$ rows 
and at most $r+1$ columns. This implies the dimension of $Z$ 
has been reduced by $r_e-r-1$. So if $r=2$, cliques larger than a 
triangle help reduce the dimension of $Z$. As mentioned above, the
existence of cliques is highly likely, since edges in the graph exist
when sensors are within radio range. Moreover, the above technique
extends to dense sets, rather than cliques. The key is finding 
 $B= (W\circ \KK)^\dagger (W\circ E_{ii})$, for an appropriate
submatrix $E_{ii}$, as well as deriving the nullspace of $W\circ \KK$.

\subsection{Estimating Sensor Positions based on \EDMt Model}
\label{sec:estimation}
After we solve $\SNL_{\scriptsize{\EDM A}}$  (or equivalently
$\SNL_{\scriptsize{\EDM s}}$)
to get an optimal solution $Z_s$, we can express
\[ \bar Y = U_s Z_s U_s^T = 
             \pmatrix{\bar Y_{11} &\bar Y_{21}^T 
                 \cr \bar Y_{21}& \bar Y_{22} }, \quad
             \bar Y_{22}= AA^T, \bar Y_{21} =  AX^T, \mbox{ for some } X.
 \]
To complete the \SNL problem, we have to find an approximation to the
matrix $P \in \MM^{n+m,r}$,      
i.e. the matrix that has the sensor locations in the first
$n$ rows, also denoted $X$, and the anchor locations in the last
$m$ rows, denoted $A$.  

Since $\bar Y \in S^{n+m}_{+}$, there exists
$\hat P = \pmatrix{\hat P_{11} & \hat P_{21} \cr \hat P_{12}  & \hat P_{22} }
\in \MM^{n+m}$ such that $\hat P \hat P^T = \bar Y$. 
By Assumption \ref{assumpt:centered}, the anchors
are centered, i.e. $A^Te=0$. We can translate the locations in
$\hat P$, so that the last $m$ locations are centered, i.e. without loss
of generality we have 
\beq
\label{eq:Pcentered}
\pmatrix{\hat P_{12}  & \hat P_{22} }^Te=0, \quad
\hat P \hat P^T = \bar Y. 
\eeq
Also
\[
\left\{\bar P\in \MM^{n+m} : \bar Y = \bar P \bar P^T\right\}= 
\left\{\bar P\in \MM^{n+m} : \bar P =\hat P Q, 
\mbox{for some orthogonal } Q \in \MM^{n+m}\right\}.
\] 
In other words, from the optimal $\bar Y$,
all the possible locations can be obtained by 
a rotation/reflection of $\hat P$. 
However, these locations in the rows of $\hat P$ are in $\R^{n+m}$, 
rather than in the
desired embedding space $\R^r$, where the anchors lie. 
\begin{rem}
Since \SNL is underdetermined, in general, the optimum $\bar Y$ 
is not unique. Therefore, finding a lower rank optimum $\bar Y$ should
result in better approximations for the sensor locations.
\end{rem}

Following are two methods for finding an
estimate to the sensor locations, $X$. The first is the one currently
used in the literature. The second is a strengthened new 
method based on the \EDM interpretation. 
\begin{enumerate}
\item
In the recent papers on \SNLc e.g. 
\cite{MR2191577,BiswasYe:04,Jin:05}, $X$ is taken directly from the optimal 
$Z_s=\pmatrix{I & X^T \cr X & Y}$,
see e.g. \eqref{eq:Zs}. Equivalently, since $A$ is full column rank $r$,
and the equations in $AX^T=\bar Y_{21}$ are consistent,
we can solve for $X$ uniquely from the $AX^T=\bar Y_{21}$.
We now describe the underlying geometry of using this $X$.

Recall that $A=U\Sigma_rV^T$ and
$\pmatrix{ \hat P_{12}  & \hat P_{22} }
 \pmatrix{ \hat P_{12}  & \hat P_{22} }^T=
  AA^T= \pmatrix{A& 0} \pmatrix{A& 0}^T$. Therefore, these three matrices
all have the same spectral decomposition and all can be 
diagonalized using $U$. This implies that the three matrices
$\pmatrix{ \hat P_{12}  & \hat P_{22} },
  A, \pmatrix{A& 0}$ can all use the same set of left singular vectors in a
compact singular value decomposition, SVD.  Therefore, 
$\pmatrix{ \hat P_{12}  & \hat P_{22} }Q= \pmatrix{A& 0}$, for some
orthogonal $Q$, i.e.
\beq
\label{eq:Xsensorsgood}
\exists \hat Q, \hat Q^T\hat Q=I, \mbox{ with }
\bar P =\hat P \hat Q  = 
\pmatrix{\bar P_{11} & \bar P_{21} \cr A  & 0 }.
\eeq
This yields
\beq
\label{eq:progrotat}
\bar P = \pmatrix{\bar P_{11} & \bar P_{12} \cr A & 0}, \quad
 \bar Y=
 \bar P \bar P^T =
\pmatrix{\bar Y_{11} & \bar P_{11}A^T \cr A\bar P_{11}^T & AA^T}.
 \eeq
Since $A\hat P_{11}^T=\bar Y_{21}=AX^T$, we see that $\hat P_{11} = X$. 
Thus the first $n$ rows of $\hat P$ project exactly onto the rows of $X$, 
after the rotation/reflection with $\hat Q$ 
to make the bottom $m$ rows equal to $A$.
If we denote the orthogonal projection onto the first $r$ coordinates by
$P_r$, then the resulting operation on the locations in the rows
of $\hat P$ can be summarized by 
\[
P_r \bar P^T=\left(P_r \hat Q^T\right) \hat P^T \in \R^r \otimes \{0\}^{n+m-r},
\bar Y \approx \bar Y_p:= \bar P P_r \bar P^T.
\]
Note that the product $P_r \hat Q^T$ is not necessarily idempotent or symmetric,
i.e. not necessarily an (orthogonal) projection. Moreover, the term that
is deleted $\bar P_{12}$ can be arbitrary large, while the rank of $\bar
Y$ can be as small as $r+1$. The relaxation from $\bar Y_{11}=XX^T$ to 
$\bar Y_{11} 
= \pmatrix{\bar P_{11} & \bar P_{12}}
 \pmatrix{\bar P_{11} & \bar P_{12}}^T=
 \bar P_{11}\bar P_{11}^T + \bar P_{12}\bar P_{12}^T
\succeq XX^T$, shows that using $X=\bar P_{11}$ has an error of the order
of $\|\bar P_{12}\|^2$.
\begin{quote} \underline{{\bf Method 1:}}
 Estimate the location of the sensors using
$X$ in the optimal $Z_s$ or, equivalently, solve for
$X$ using  the equation $AX^T=\bar Y_{21}$, where 
$\bar Y_{21}$ is from the optimal $\bar Y$.
\end{quote}

\item
In Method 1, the matrix $\bar P P_r \bar P^T$ provides 
a $\rank r$ approximation to $\bar Y$.
However, if $\|\bar P_{12}\|$ in \eqref{eq:progrotat} is large, 
then it appears that we have lost information.
It is desirable to keep as much of the information 
from the high dimensional locations in $\hat P$ as we can, i.e. the
information that is contained in $\bar Y_{11}$.
If we do not consider the anchors distinct from the sensors,
then we would like to rotate and then project
all the rows of $\hat P$ onto a subspace of
dimension $r$, i.e. we consider the problem to be an \EDM completion
problem and would like to extract a good approximation of the positions
of all the nodes.
 Since the last $m$ rows corresponding to the anchors
originated from a clique, the corresponding graph is rigid and the
corresponding projected points will be close to the original anchor
positions.
We realize this using the spectral decomposition.
(See e.g. \cite{AlKaWo:97}, where error estimates are included.)
\[
\bar Y = \pmatrix{U_r &\bar U_r} \pmatrix{\Sigma_r & 0 \cr 0 & 
           \Sigma_{n+m-r}}\pmatrix{U_r & \bar U_{r}}^T.
\]
Then, considering the problem as an \EDM completion problem, we {\em
first} find a {\em best} rank $r$ approximation to $\bar Y$,  
denoted $\bar Y_r:=U_r \Sigma_r U_r^T$. Only then do
we find a particular full rank factorization $\hat P_r \in \MM^{n+m,r}$ 
such that $\bar Y_r =\hat P_r \hat P_r^T$, i.e.
$\hat P_r = U_r \Sigma_r^{1/2}$.  It remains to find an
orthogonal $Q$ in order to find $P=\hat P_r Q$.
Fortunately, we can use the information from the anchors to find the
orthogonal $Q$. 
\begin{quote} \underline{{\bf Method 2:}}
Suppose that $\hat P_r = U_r \Sigma_r^{1/2}$ is found as above, with
$\hat P_r = \pmatrix{\hat P_1 \cr \hat P_2}$.
We find $\hat Q$ as a minimum for
$\min_{Q^TQ=I} \|\hat P_2 Q - A \|^2_F$. 
The solution is given analytically by
$\hat Q = V_QU_Q^T$, where $U_Q\Sigma_Q V_Q=A^T \hat P_2$ is the 
SVD for $A^T \hat P_2$. 
Then the rows of 
$\hat P_1 \hat Q$ are used to estimate the locations of the sensors.
\end{quote}
\end{enumerate}

Numerical tests for the two methods, are given in 
Section  \ref{sect:numerestsensors}.
Method 2 proved to be consistently more accurate.
However, method 1 locates all sets of sensors that are uniquely 
localizable in $R^r$, see \cite{SoYe:05}.

\begin{rem}
As above, suppose that $\bar Y$ is an optimum for the \SDP relaxation.
The problem of finding a best $P$ to estimate the sensor locations 
is equivalent to finding
\[ 
P^* \in \argmin_{P=\pmatrix{X \cr A}} 
 \| W \circ (\KK(PP^T-\bar Y))\|_F.
\]
Equivalently, we want to find
\[
Y^* \in S^{n+m}_{+}, \rank (Y^*)= r, Y^*_{22} = AA^T, Y^* = 
        \bar Y+\NN(W \circ \KK).
\]
However, finding such a $Y^*$ is equivalent to finding the minimal rank 
matrix in the intersection of the semidefinite cone and an affine space. 
This is still an open/hard problem. Recently, 
\cite{SoYeZhang:06,LuoSidTsengZhang:06} 
proposed randomization methods for \SDP rank reduction. 
These methods can generate a low rank positive semidefinite 
matrix in an approximate affine space. 
\end{rem}

\section{Duality for \SNLt with Quadratic Constraint}
\label{sect:dualnonlin} 
Instead of using the standard linearized relaxation as 
\SNLMV in \eqref{eq:relax1} and in \cite{KrPiWo:06}, 
we now study the new relaxation without linearizing the quadratic constraint 
$XX^T - Y \preceq 0$. This avoids ill-conditioning caused by this
linearization, see Remark \ref{rem:onto}.
Our numerical results indicate that the new quadratic approach is 
more stable than the linear approach, 
see more in Section \ref{sect:numeric}.
A discussion on the strengths of the corresponding barriers is given in
\cite{GulerTuncel:98,ChuaTu:05}.

Recall that $x=\sqrt{2}\kvec(X),~y:=\svec(Y)$. We begin with the reduced 
problem
\beq \label{eq:relaxNonlinear} \mbox{(\SNLMN)} \qquad
  \begin{array}{crcl}
      \min  &  f_3(x,y) :=  \frac{1}{2}{\| W \circ (\KK(\YY(x,y))) - \bar E
\|}_F^2 \\
  \mbox{subject to } & g_u(x,y):=H_u \circ \KK(\YY(x,y)) - \bar U &\leq& 0
\\
                   & g_l(x,y):=\bar L -H_l \circ \KK(\YY(x,y))& \leq& 0 \\
                & \frac{1}{2}\Mat(x){\Mat(x)}^T - \sMat(y) \preceq 0.
  \end{array}
\eeq\\
Then the Lagrangian is
\beq \label{eq:lagrangianNonlinear}
\begin{array}{rcl}
    L(x,y,\Lambda_u,\Lambda_l,\Lambda)
    &=&
      \frac{1}{2}{\| W \circ \KK(\YY(x,y)) - \bar E \|}_F^2
         + \ip{\Lambda_u}{H_u \circ \KK(\YY(x,y))-\bar U}
  \\&& \quad
            +  \ip{\Lambda_l}{\bar L -H_l \circ \KK(\YY(x,y))}
  \\&& \quad
       + \ip{\Lambda}{\frac{1}{2}\Mat(x){\Mat(x)}^T - \sMat(y)},
\end{array}
\eeq 
where $0\leq \Lambda_u, 0\leq \Lambda_l \in \Ss^{m+n}$, and
$0\preceq \Lambda \in \Ss^{n}$. In addition, we denote
\[
\begin{array}{c}
                \lambda_u:=\svec(\Lambda_u), \quad
                      \lambda_l:=\svec(\Lambda_l), \\
                h_u:=\svec(H_u), \quad
                      h_l:=\svec(H_l),\quad
\lambda:=\svec(\Lambda).
\end{array}
\]
And, for numerical implementation, we define the linear
transformations \beq \label{eq:hnz2} h_u^{nz}= \svec_u ( H_u) \in
\R^{nz_u}, \quad h_l^{nz}= \svec_l ( H_l) \in \R^{nz_l}, \eeq where
$h_u^{nz}$ is obtained from $h_u$ by removing the zeros; thus,
$nz_u$ is the number of nonzeros in the upper-triangular part of
$H_u$. Thus the indices are fixed from the given matrix $H_u$.
Similarly, for $h_l^{nz}$ with indices fixed from $H_l$. We then get
the vectors
\[
\lambda_u^{nz}= \svec_u ( \Lambda_u) \in \R^{nz_u}, \quad
\lambda_l^{nz}= \svec_l ( \Lambda_l) \in \R^{nz_l}.
\]
The adjoints are $\sMat_u,\sMat_l$; and, for any matrix $M$ we get
\[
H_u\circ M= \sMat_u \svec_u ( H_u\circ M).
\]
This holds similarly for $H_l\circ M$. Therefore, we could rewrite
the Lagrangian as \beq \label{eq:lagrangianbig2b}
\begin{array}{rcl}
    L(x,y,\Lambda_u,\Lambda_l,\Lambda)
    &=&
    L(x,y,\lambda^{nz}_u,\lambda^{nz}_l,\Lambda)
  \\&=&
      \frac{1}{2}{\| W \circ \KK(\YY(x,y)) - \bar E \|}_F^2
         + \ip{\svec_u(\Lambda_u)}
                {\svec_u\left(H_u \circ \KK(\YY(x,y))-\bar U\right)}
  \\&& \quad
            +  \ip{\svec_l(\Lambda_l)}{\svec_l\left(
          \bar L -H_l \circ \KK(\YY(x,y))\right)}
  \\&& \quad
       + \ip{\Lambda}{\frac{1}{2}\Mat(x){\Mat(x)}^T - \sMat(y)}.
\end{array}
\eeq\\
To simplify the dual of \SNLMN, i.e. the max-min of the Lagrangian,
we now find the stationarity conditions of the inner minimization
problem, i.e. we take the derivatives of $L$ with respect to $x$ and
$y$. We get
\beq
\label{eq:derivXnonlin}
\begin{array}{rcl}
0 &=&
    \nabla_x L(x,y,\Lambda_u,\Lambda_l,\Lambda)
     \\&=&
  \left[W\circ(\KK \YY^x)\right]^*
       \left( W \circ \KK(\YY(x,y)) - \bar E \right)
     +
  \left[H_u\circ(\KK \YY^x)\right]^* (\Lambda_u)
  \\&& \hspace{2.5in}  -
  \left[H_l\circ(\KK \YY^x)\right]^* (\Lambda_l)
        +
       \kvec(\Lambda \Mat(x)).
\end{array}
\eeq
Note that
\beq
\begin{array}{rcl}
  T(x) &=&  \ip{\Lambda}{\frac{1}{2}\Mat(x){\Mat(x)}^T} \\
  &=& \frac{1}{2}\ip{x}{\kvec(\Lambda\Mat(x))}.
\end{array}
\eeq Therefore,  $\frac{dT(x)}{dx} = \kvec(\Lambda\Mat(x))$, since
$\ip{x}{\kvec(\Lambda\Mat(x))}$ is a quadratic form in x.
Similarly,
\beq\label{eq:derivYnonlin}
\begin{array}{rcl}
0 &=&
    \nabla_y L(x,y,\Lambda_u,\Lambda_l,\Lambda)
     \\&=&
  \left[W\circ(\KK \YY^y)\right]^*
       \left( W \circ \KK(\YY(x,y)) - \bar E \right)
     +
  \left[H_u\circ(\KK \YY^y)\right]^* (\Lambda_u)
  \\&& \hspace{2.5in}  -
  \left[H_l\circ(\KK \YY^y)\right]^* (\Lambda_l)
        -
       \svec(\Lambda),
\end{array}
\eeq
since $\ip{\Lambda}{\sMat(y)}$ is linear in y.
We can solve for $\Lambda$ and then use this to eliminate $\Lambda$ in the
other optimality conditions, i.e. we eliminate $t(n)$ variables and
equations using
\beq
\label{eq:derivYnonlinLambda}
\begin{array}{rcl}
       \svec(\Lambda)
&=&
  \left[W\circ(\KK \YY^y)\right]^*
       \left( W \circ \KK(\YY(x,y)) - \bar E \right)
     +
  \left[H_u\circ(\KK \YY^y)\right]^* (\Lambda_u)
  \\&& \hspace{2.5in}  -
  \left[H_l\circ(\KK \YY^y)\right]^* (\Lambda_l).
\end{array}
\eeq
We now substitute for $\Lambda$ in the first stationarity condition
\eqref{eq:derivXnonlin}, i.e.
\beq
\label{eq:derivXnonlinelim}
\begin{array}{rcl}
0 &=&
  \left[W\circ(\KK \YY^x)\right]^*
       \left( W \circ \KK(\YY(x,y)) - \bar E \right)
     +
  \left[H_u\circ(\KK \YY^x)\right]^* (\Lambda_u)
  -\left[H_l\circ(\KK \YY^x)\right]^* (\Lambda_l)
  \\&& \hspace{.1in}
        + \kvec\left( \sMat\left\{
  \left[W\circ(\KK \YY^y)\right]^*
       \left( W \circ \KK(\YY(x,y)) - \bar E \right) \right. \right.
   \\&& \hspace{.1in}  
     +
  \left. \left. \left[H_u\circ(\KK \YY^y)\right]^* (\Lambda_u)
      -
  \left[H_l\circ(\KK \YY^y)\right]^* (\Lambda_l)\right\}
\Mat(x)\right).
\end{array}
\eeq

The Wolfe dual is obtained from applying the stationarity conditions
to the inner minimization of the Lagrangian dual (max-min
of the Lagrangian), i.e. we get the (dual \SNLMN) problem
\beq
\label{eq:dualrelaxNonlinear} (\SNLMND)
\begin{array}{rcl}
      &\mbox{max} &
  L(x,y,\lambda_u,\lambda_l,\lambda)\\
       &\mbox{subject to }&
\eqref{eq:derivXnonlin},\eqref{eq:derivYnonlin}\\
        && \sMat(\lambda_u) \geq 0, \sMat(\lambda_l) \geq 0 \\
                          &&
         \sMat (\lambda) \succeq 0.
  \end{array}
\eeq
We denote the slack variables \beq \label{eq:primfeas2b}
\begin{array}{rcl}
S_u & := & \bar U - H_u \circ \left(
\KK\left(\YY(x,y)\right)\right),\quad
          s_u=\svec S_u \\
  S_l & := & H_l \circ \left(\KK\left(\YY(x,y)\right)\right) -\bar L,\quad
          s_l=\svec S_l\\
Z & := & Y - XX^T \succeq 0.
\end{array}
\eeq \\
We can now present the primal-dual characterization of
optimality.
\begin{thm}
\label{thm:optcond12} The primal-dual variables
$x,y,\Lambda,\lambda_u,\lambda_l$ are optimal for \SNLMN if and only if:
\begin{enumerate}
\item Primal Feasibility:
\[
s_u \geq 0, \quad s_l \geq 0, \mbox{ in \eqref{eq:primfeas2b}},
\]
\beq \label{eq:pdfeas2}
    \frac 12\Mat(x){\Mat(x)}^T - \sMat(y) \preceq 0.
\eeq
\item Dual Feasibility: Stationarity equations
\eqref{eq:derivXnonlin},\eqref{eq:derivYnonlin} hold and \beq
\label{eq:sdpL1}
  \Lambda=
     \sMat(\lambda) \succeq 0; \lambda_u\geq 0;\lambda_l\geq 0. \eeq
\item Complementary Slackness:
\beq
\label{eq:complslackF}
\begin{array}{rcl}
    \lambda_u \circ s_u & = & 0 \\
    \lambda_l \circ s_l & = & 0 \\
    \Lambda Z & = & 0.
\end{array}
\eeq

\end{enumerate}
\epr
\end{thm}
We can use the structure of the optimality conditions  to eliminate
some of the linear dual equations and obtain a characterization of
optimality based mainly on a bilinear equation and
nonnegativity/semidefiniteness.
\begin{cor}
The dual linear equality constraints \eqref{eq:derivYnonlin} in
Theorem \ref{thm:optcond12}  can be eliminated after using it to
substitute for $\lambda$ in \eqref{eq:derivXnonlin}, i.e. we get
equation \eqref{eq:derivXnonlinelim}. The complementarity conditions
in \eqref{eq:complslackF} now yield a bilinear system of equations
$F(x,y,\lambda_u,\lambda_l)=0$,
with nonnegativity and semidefinite conditions that characterize
optimality of \SNLMN.
\epr
\end{cor}

\section{A Robust Primal-Dual Interior-Point Method}
\label{sect:pdipNonlin} 
We now present a primal-dual interior-point method for \SNLMN, see in \cite{KrPiWo:06} for the linearized case, \SNLMV.
First, we define the equation
\eqref{eq:derivXnonlinelim} to be:
\[
L_s(x,y,\Lambda_u,\Lambda_l) = 0.
\]
Then, to solve \SNLMN we use the Gauss-Newton method on the
perturbed complementary slackness conditions (written with the block
vector notation): \beq \label{eq:pertcsNonlin}
\begin{array}{rcl}
F_\mu(x,y,\lambda_u,\lambda_l):= \pmatrix{
    \lambda_u \circ s_u  - \mu_u e\cr
    \lambda_l \circ s_l  - \mu_l e \cr
    \Lambda Z  -  \mu_c I \cr
    L_s } = 0,
\end{array}
\eeq where $s_u=s_u(x,y)$, $s_l=s_l(x,y)$,
$\Lambda=\Lambda(x,y,\lambda_u,\lambda_l)$, $Z=Z(x,y)$ and $L_s =
L_s(x,y,\Lambda_u,\Lambda_l).$ This is an overdetermined system with
\[
(m_u+n_u)+ (m_l+n_l)+ n^2 + nr \mbox{ equations}; \quad nr +t(n) +
(m_u+n_u)+ (m_l+n_l) \mbox{ variables}.
\]

\subsection{Linearization}
We denote the Gauss-Newton search direction for
\eqref{eq:pertcsNonlin} by
\[
    \Delta s := \left(\begin{array}{c}
        \Delta x \\\Delta y \\  \Delta \lambda_u \\ \Delta \lambda_l
    \end{array}\right).
\]

The linearized system for the search direction $\Delta s$ is:
\[
F^\prime_\mu(\Delta s) \cong F^\prime_\mu(x,y,\lambda_u,\lambda_l)
(\Delta s) =-F_\mu(x,y,\lambda_u,\lambda_l).
\]

To further simplify notation, we use the following composition of
linear transformations. Let $H$ be symmetric. Then
\[
\begin{array}{rcl}
\KK^x_H (x)&:=& H\circ (\KK (\YY^x(x))),\\
\KK^y_H (y)&:=& H\circ (\KK (\YY^y(y))),\\
\KK_H (x,y)&:=& H\circ (\KK (\YY(x,y))).
\end{array}
\]
so, we have the following:
\[
\begin{array}{rcl}
\Lambda (x,y,\lambda_u,\lambda_l) & = & \sMat [(\KK^y_W)^*(K_W(x,y)-
\bar E) + (\KK^y_{H_u})^*(\sMat(\lambda_u)) -
(\KK^y_{H_l})^*(\sMat(\lambda_l))],\\
L_s(x,y,\Lambda_u,\Lambda_l) & = & (\KK^x_W)^*(K_W(x,y)- \bar E) +
(\KK^x_{H_u})^*(\sMat(\lambda_u)) -
(\KK^x_{H_l})^*(\sMat(\lambda_l)) \\
&& \quad + \kvec (\Lambda \Mat (x)).
\end{array}
\]
Define the linearization of above functions as:
\[
\begin{array}{rcl}
\Delta \Lambda (\Delta x,\Delta y,\Delta \lambda_u,\Delta \lambda_l)
& = & \sMat [(\KK^y_W)^*(K_W(\Delta x, \Delta y)) +
(\KK^y_{H_u})^*(\sMat(\Delta \lambda_u))\\
&&\quad - (\KK^y_{H_l})^*(\sMat(\Delta \lambda_l))],\\
\Delta L_s(\Delta x,\Delta y,\Delta \Lambda_u,\Delta \Lambda_l) &=&
(\KK^x_W)^*(K_W(\Delta x,\Delta y)) + (\KK^x_{H_u})^*(\sMat(\Delta
\lambda_u)) - (\KK^x_{H_l})^*(\sMat(\Delta \lambda_l)) \\
&& \quad + \kvec (\Delta \Lambda \Mat (x)) + \kvec (\Lambda \Mat
(\Delta x)).
\end{array}
\]
The linearization of the complementary slackness conditions results
in four blocks of equations
\begin{enumerate}
\item
\[
- \lambda_u \circ \svec \KK_{H_u} (\Delta x,\Delta y) + s_u \circ
\Delta \lambda_u = \mu_u e - \lambda_u \circ s_u
\]
\item
\[
\lambda_l \circ \svec \KK_{H_l} (\Delta x,\Delta y) + s_l \circ
\Delta \lambda_l = \mu_l e - \lambda_l \circ s_l
\]
\item
\begin{eqnarray*}
&&\Lambda (\sMat(\Delta y) - \frac{1}{2}\Mat(x)\Mat(\Delta x)^T -
\frac{1}{2}\Mat(\Delta x)\Mat(x)^T)\\
&& + \Delta \Lambda (\Delta s) (\sMat(y) - \frac{1}{2}\Mat(x)\Mat(x)^T) \\
&=& \mu_cI-\Lambda Z
\end{eqnarray*}
\item
\[
\Delta L_s(\Delta s) = -L_s(x,y,\Lambda_u,\Lambda_l)
\]
\end{enumerate}
and hence
\[
    F^\prime_\mu (\Delta s)  = \\
    \left(
    \begin{array}{l}
        - \lambda_u \circ \svec \KK_{H_u} (\Delta x,\Delta y) + s_u \circ
\Delta \lambda_u \\
        \lambda_l \circ \svec \KK_{H_l} (\Delta x,\Delta y) + s_l \circ
\Delta \lambda_l \\
        \Lambda (\sMat(\Delta y) - \frac{1}{2}\Mat(x)\Mat(\Delta x)^T - 
\frac{1}{2}\Mat(\Delta x)\Mat(x)^T)\\
\quad + \Delta \Lambda (\Delta s)Z \\
        \Delta L_s(\Delta s)
    \end{array}
    \right)
\]
where $F_\mu^\prime:\MM^{n \times r}\times\Re^{t(n)}\times
\Re^{t(m+n)}\times \Re^{t(m+n)}\rightarrow
\Re^{t(m+n)}\times\Re^{t(m+n)}\times\MM^{n \times n}\times\Re^{nr}$,
i.e. the linear system is overdetermined.

We need to calculate the adjoint $(F^\prime_\mu)^*$. We first find
$(\KK_H^x)^*$, $(\KK_H^y)^*$, and $(\KK_H)^*$. By the expression of
$\YY(\Delta x,\Delta y)$, we get \beq \YY^*(S)
=\pmatrix{(\YY^x)^*(S)\cr (\YY^y)^*(S)}=
\pmatrix{\Mat^*(\sblk_{21}(S)^TA)\cr\sMat^*(\sblk_1(S)) }=\pmatrix{
\kvec(\sblk_{21}(S)^TA)\cr\svec(\sblk_{1}(S))}. \eeq By the
expression of $\KK_H(\Delta x,\Delta y)$, we get \beq \KK_H^*(S)
=\pmatrix{(\KK_H^x)^*(S)\cr (\KK_H^y)^*(S)}=
\pmatrix{(\YY^x)^*(\KK^*(H\circ S))\cr(\YY^y)^*(\KK^*(H\circ S)) }.
\eeq 
Moreover,
\begin{eqnarray*}
\ip{\Lambda \sMat (\Delta y)}{W_3}
&=&
\tr \left(W_3^T \Lambda\right)\sMat (\Delta y)
\\ &=&\ip{\frac{1}{2}\svec (\Lambda W_3 + W_3^T \Lambda)}{\Delta y}
\end{eqnarray*}
Similarly
\begin{eqnarray*}
\ip{\frac{1}{2} \Lambda \Mat(x) \Mat(\Delta x)^T}{W_3}
&=&
\tr \frac{1}{2}W_3^T \Lambda \Mat(x) {\Mat(\Delta x)}^T 
\\&=& \ip{\frac{1}{2}\kvec (W_3^T \Lambda \Mat(x))}{\Delta x}
\end{eqnarray*}
and
\begin{eqnarray*}
\ip{\frac{1}{2} \Lambda \Mat(\Delta x) \Mat(x)^T}{W_3}
&=&
\tr \frac{1}{2}W_3^T \Lambda \Mat(\Delta x){\Mat(x)}^T 
\\&=& \ip{\frac{1}{2}\kvec (\Lambda W_3 \Mat(x))}{\Delta x}.
\end{eqnarray*}
We also need to find $(\Delta \Lambda)^*(S)$ by the expression of
$\KK_H(\Delta x,\Delta y)$, where $S\in\Sn $, we get \beq (\Delta
\Lambda)^*(S) = \pmatrix{(\KK_W^x)^*(\KK_W^y(\svec (S)))\cr
(\KK_W^y)^*(\KK_W^y(\svec (S)))\cr \svec[(\KK_{H_u}^y)(\svec
(S))]\cr -\svec[(\KK_{H_l}^y)(\svec (S))]} = \pmatrix{ \Delta x\cr
\Delta y\cr \Delta \lambda_u\cr \Delta \lambda_l}. 
\eeq 
Then we have
\begin{eqnarray*}
\ip{\Delta \Lambda (\Delta x,\Delta y,\Delta \lambda_u,\Delta 
\lambda_l)(Z)}{W_3}
&=& \tr W_3^T \Delta \Lambda (\Delta s) (Z)
\\&=& \ip{(\Delta \Lambda)^*(\frac {1}{2} [W_3Z+ZW_3^T])}{\pmatrix{\Delta 
x\cr\Delta y\cr\Delta \lambda_u\cr\Delta \lambda_l}}.
\end{eqnarray*}
Now we find $(\Delta L_s)^*(w_4)$, which consists of three columns
of block with four rows per column. We list this by columns $C_1,
C_2,C_3$. 
\begin{eqnarray*} C_1 =
\pmatrix{(\KK_W^x)^*(\KK_W^x(w_4))\cr (\KK_W^y)^*(\KK_W^x(w_4))\cr
\svec[(\KK_{H_u}^x)(w_4)]\cr -\svec[(\KK_{H_l}^x)(w_4)]}
\end{eqnarray*}

\begin{eqnarray*} C_2 = \pmatrix{\kvec(\Lambda \Mat(w_4))\cr 0 \cr 0 \cr 0 
}\end{eqnarray*}

\begin{eqnarray*} C_3 = (\Delta \Lambda)^*(\frac{1}{2}[\Mat(w_4)X^T + 
X\Mat(w_4)^T])\end{eqnarray*}
Thus, the desired adjoint is given by $(\Delta L_s)^*(w_4) = C_1 +
C_2 + C_3$.

Now we evaluate $(F^\prime_\mu)^*(w_1,w_2,W_3,w_4)$. This consists
of four columns of blocks with four rows per column. We list this by
columns $Col_1, Col_2, Col_3, Col_4$.
\begin{eqnarray*}
Col_1 &=&{\pmatrix{
        -(\KK^x_{H_u})^*\left(\sMat(\lambda_u\circ
w_1)\right)\cr
    -(\KK^y_{H_u})^*\left(\sMat(\lambda_u\circ
w_1)\right)
                       \cr    
         w_1\circ s_u
         \cr                 
            0
         }}
\end{eqnarray*}
\begin{eqnarray*}
Col_2&=&{\pmatrix{
       (\KK^x_{H_l})^*\left(\sMat(\lambda_l\circ
w_2)\right)\cr
    (\KK^y_{H_l})^*\left(\sMat(\lambda_l\circ
w_2)\right)
                       \cr    
         0\cr                 
         w_2\circ s_l\cr       
}}
\end{eqnarray*}
\begin{eqnarray*} Col_3 = Col_{31}+Col_{32}\end{eqnarray*}
\begin{eqnarray*} Col_4 = (\Delta L_s)^*(w_4)\end{eqnarray*}
where
\begin{eqnarray*}
Col_{31} &=& {\pmatrix{-\frac{1}{2}\kvec 
   \left(W_3^T \Lambda \Mat(x)+
            \Lambda W_3 \Mat(x)\right)\cr \frac{1}{2}\svec (\Lambda
W_3 + W_3^T \Lambda)\cr 0 \cr 0 } }
\end{eqnarray*}
\begin{eqnarray*} Col_{32} = (\Delta \Lambda)^*(\frac{1}{2}[W_3Z+ZW_3^T 
])\end{eqnarray*}
where $w_1\in\Re^{t(m+n)},w_2\in\Re^{t(m+n)},W_3\in\MM^{n \times n}$
and $w_4\in\Re^{nr}$. Thus the desired adjoint is given by
$(F_\mu^\prime)^*=Col_1+Col_2+Col_3+Col4$.

\section{Numerical Tests}
\label{sect:numeric}
We now present results on randomly generated 
\SNL problems with connected underlying graphs. The tests were done
using MATLAB 7.1.
The method for generating the tests follows from the approach used in 
\cite{Jin:05,KrPiWo:06}.

The first set of tests compares the two methods for finding a proper
factorization to estimate the sensor locations from the optimum of the
\SDP relaxation. The second set of tests compares the two methods 
for solving the \SDP relaxation, i.e. using the quadratic constraint $XX^T
\preceq Y$  and the linear one using $Z_s\succeq 0$.

\subsection{Two Methods for Estimating Sensor Locations}
\label{sect:numerestsensors}
Two methods for estimating the sensor locations from a given optimum of
the \SDP relaxation were presented in Section \ref{sec:estimation}, i.e.
\begin{enumerate}
\item{Method 1:}
After obtaining $Z_s = \pmatrix{I_r & X^T \cr X & Y}$, $X$ is used 
to estimate sensor positions as in \cite{MR2191577,BiswasYe:04,Jin:05}.
\item{Method 2:} 
Use the rows of 
$\hat P_1 \hat Q$ to estimate the locations of the sensors. Here
$\hat Q = V_QU_Q^T$, $U_Q\Sigma_Q V_Q=A^T \hat P_2$ is the 
singular value decomposition for $A^T \hat P_2$. 
\end{enumerate}

We denote $X_{e1},X_{e2}$ as the estimated sensor locations from method 1 
and method 2, respectively. 
We first note there is a significant difference in norm between the estimates
of $X$ from Method 1 and Method 2, see Table \ref{table:comptest4}.
In Tables
\ref{table:comptest1},\ref{table:comptest2},\ref{table:comptest3},
we use the following three measures to compare the two methods.
\begin{description}
\item {Measure 1: Objective Function with Different Anchors}\\
When finding the SVD decomposition of $AA^T$
in method 2, the anchor locations are estimates, $A_{e2}$, 
and may not correspond to $A$.
This measure uses the true objective function 
$\left\|W \circ \left(\KK\left(\pmatrix{X_{ei}X_{ei}^T & X_{ei}A_{ei}^T 
     \cr A_{ei}X_{ei}^T & A_{ei}A_{ei}^T}\right)-E\right)\right\|_F, i=1,2$.

\item {Measure 2: Total Distance Error with True Sensor Locations}\\
This is a common criterion used for \SNL. 
We compare the sum of distances between estimated sensor 
locations and true sensor locations, i.e. $ \|X_{ei}-X^*\|_F,i=1,2$, 
where $X^*$ denotes the true sensor locations. 

\item {Measure 3: Objective Function with Original Anchors}\\
We use the same criterion as in Measure 1, except that we keep the
anchors  fixed to their true locations, i.e.
$\left\|W \circ \left(\KK\left(\pmatrix{X_{ei}X_{ei}^T & X_{ei}A^T 
     \cr AX_{ei}^T & AA^T}\right)-E\right)\right\|_F, i=1,2$.
\end{description}

\begin{table}[htbp]
\centering
{\small
\begin{tabular}{|c|c|c|c|c|c|c|c|c|c|}
\hline
 &test 1&test 2 &test 3&test 4 &test 5&test 6&test 7&mean& std  \\
\hline
$\|X_{e1}-X_{e2}\|_F$ & 1.2351 &1.3002 & 1.4210 &1.2906&1.1300&1.3810&1.2964&1.2935 &0.0879\\
\hline
\end{tabular}
}
\caption{
Difference of the estimates from the two methods}
\label{table:comptest4}
\end{table}~\begin{table}[htbp] 
\centering
{\small
\begin{tabular}{|c|c|c|c|c|c|c|c|c|c|}
\hline
 &test 1&test 2 &test 3&test 4 &test 5&test 6& test 7&mean & std \\
\hline
Method 1 & 3.5200&3.8549&3.8793 & 3.5006 & 2.9434&3.4693&3.8736&3.5773&0.3363  \\
\hline
Method 2 &0.7462  &0.9213&0.8656 & 1.0310 & 0.7237&1.6671&1.2351&1.0271 & 0.3319 \\
\hline
\end{tabular}
}
\caption{
Measure 1: Use $X$ estimates in objective function; with different anchors
}
\label{table:comptest1}
\end{table}~\begin{table}[htbp] 
\centering
{\small
\begin{tabular}{|c|c|c|c|c|c|c|c|c|c|}
\hline
 &test 1&test 2 &test 3&test 4 &test 5 & test 6 &test 7 &mean & std\\
\hline
Method 1 & 1.2780 &1.4200 & 1.4801 &1.3696&1.1820&1.4317&1.3912&1.3647& 0.1021\\[3pt]
\hline
Method 2 &0.1887  &0.1630 &0.1050 &0.1394 & 0.0778&0.0808&0.3881&0.1633&0.1074\\[3pt]
\hline
\end{tabular}
}
\caption{
Measure 2: Use distance of $X$ estimates from true sensor locations $X^*$
}
\label{table:comptest2}
\end{table}~\begin{table}[htbp] 
\centering
{\small
\begin{tabular}{|c|c|c|c|c|c|c|c|c|c|}
\hline
 &test 1&test 2 &test 3&test 4 &test 5&test 6&test 7&mean &std  \\
\hline
Method 1 & 3.5200 &3.8549 & 3.8793 &3.5006&2.9434&3.4693&3.8736&3.5773&0.3363 \\[3pt]
\hline
Method 2 &0.2771  &0.3264 &0.1588 &0.1799 & 0.1714&0.1453&0.6428&0.2717&0.1770\\[3pt]
\hline
\end{tabular}
}
\caption{
Measure 3: Use $X$ estimates in objective function; with original anchors
}
\label{table:comptest3}
\end{table}

In the tests in 
Tables \ref{table:comptest1},\ref{table:comptest2},\ref{table:comptest3},
we used randomly generated graphs with parameters:
$r=2, n=16, m=5$, and radio range $0.15$. The density of edges that are
known was $0.75$ and all the sensors/anchors lie within a $2 \times 2$ 
square. 
We also tested many instances with different parameters, e.g. more 
sensors, and larger radio range. But the results of comparing
the two methods were similar to those presented in these tables.

\subsection{Two Methods for Solving \SNLt}
\label{sect:numerbarriers}
\begin{figure}[htb]
\epsfxsize=350pt
\centerline{\epsfbox{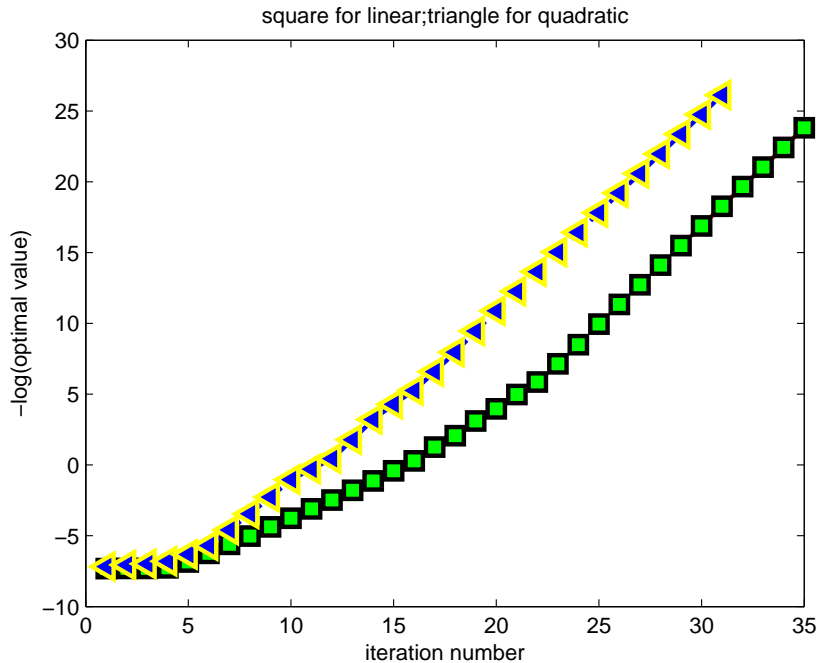}}
\caption{Comparison for two barriers; optimal value}
\label{fig:logOptval}
\end{figure}

\begin{figure}[htb]
\epsfxsize=350pt
\centerline{\epsfbox{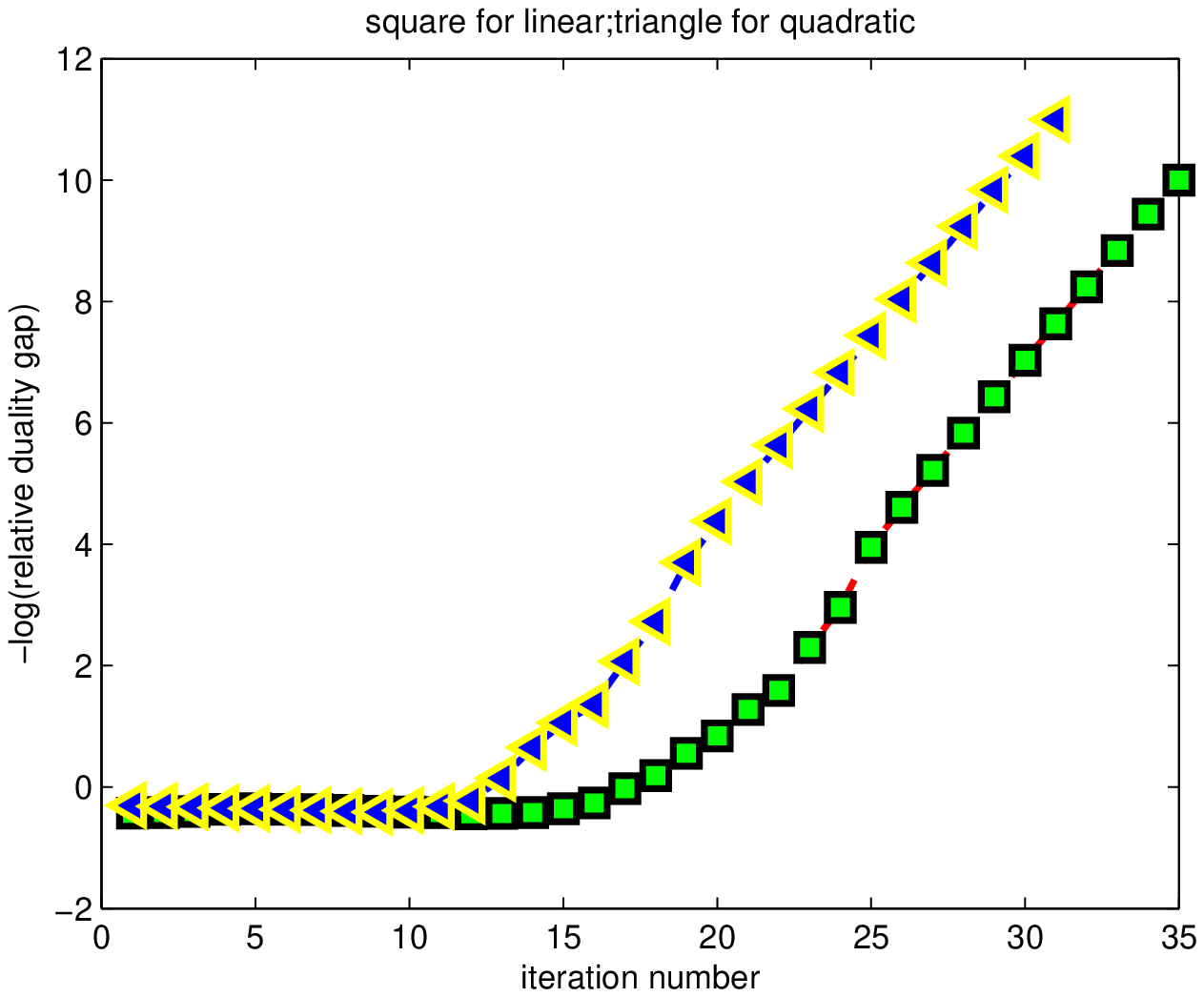}}
\caption{Comparison for two barriers; relative gap}
\label{fig:relGap}
\end{figure}

In Figures \ref{fig:logOptval},\ref{fig:relGap}, 
we present results for using the quadratic
constraint $Y-XX^T \succeq 0$ compared to the linearized version
$\pmatrix{I & X^T\cr X & Y} \succeq 0$. We solved many 
randomly generated problems with
various values for the parameters. We present typical results in the
figures.

Figure \ref{fig:logOptval} shows the ($-log$) of the optimal
value at each iteration. 
Figure \ref{fig:relGap} shows the ($-log$) of the relative gap.
Both figures illustrate the surprising result that the quadratic
formulation is more efficient, i.e. it obtains higher accuracy with fewer
iterations.  This is surprising, since we are using a Newton based method
that should be faster on functions that are less nonlinear.
Therefore, from a numerical analysis viewpoint, it appears that the
linear version is more ill-conditioned, as was mentioned since the
constraint is not {\em onto}. In addition, the figures show the high
accuracy that can be obtained though these problems are highly
ill-conditioned.

These tests provide empirical evidence for the 
theoretical comparison results
on different barriers given in \cite{GulerTuncel:98,ChuaTu:05}. 
The results in these references show that the central path is distorted
due to the $I$ in the linear formulation constraint. And, the distortion
increases with increasing dimension of the $I$.
This agrees with our
interpretation that the linear constraint is not {\em onto}, and the
Jacobian is singular.

\section{Concluding Remarks}
\label{sect:concl}
In this paper, we have analyzed the well known \SNL problem from a 
new perspective. By considering the set of anchors as a
clique in the underlying graph,  the \SNL problem can be studied using
traditional \EDM theory. Our main contributions follow from this
\EDM approach:
\begin{enumerate}
\item 
The Slater constraint qualification can fail for  cliques and/or dense
subgraphs in the underlying graph. If this happens, then
we can project the feasible set of the \SDP relaxation 
to the {\em minimal cone}. This projection improves the stability and 
can also reduce the size of the \SDP significantly. 

In a future study we plan on identifying the appropriate dense subgraphs.
(Algorithms for finding dense subgraphs exist in the literature, e.g.
\cite{MR2140670,MR2012813,MR1677349}.)

\item We provided a geometric interpretation for the method of 
directly using the $X$ from the optimal $Z_s$ of the \SDP relaxation,
when estimating the sensor positions. We then proposed another method of 
estimating the sensor positions based on a principal component analysis.
Our numerical tests showed that the new method 
gave consistently more accurate solutions. 
\item We used the $\ell_2$ norm formulation instead of the $\ell_1$ norm. This
is a better fit for the data that we used. However, the quadratic
objective makes the problem more difficult to solve.

In the future we plan on completing an error analysis comparing the two
norms.
\item
We solved the $\ell_2$ norm formulation of the \SDP relaxation
with a Gauss-Newton primal-dual interior-exterior path following
method. This was a robust approach 
compared with the traditional symmetrization and a Newton method. 
We compared using the quadratic constraint with the linearized version
used in the literature.
The numerical results showed that the quadratic constraint
is more stable. This agrees with theoretical results 
in the literature on the deformation of the central path based on the
size of the $I$ in the linearized version.

Future work involves making the algorithm more efficient. In particular,
this requires finding appropriate preconditioners.
\end{enumerate}

\bibliography{.master,.psd,.edm,.bjorBOOK}

\end{document}